\documentclass[12pt]{article}

\usepackage{tikz}
\usepackage{color}
\usepackage{amsmath}
\usepackage{amssymb}
\usepackage{enumitem}
\usepackage{graphicx}

\def\vir{\mathrm {vir}}
\def\ct{\mathrm{ct}}
\def\P{\mathsf{P}}
\def\PP{\mathbb{P}}
\def\DR{\mathsf{DR}}
\def\mathsff{{\mathcal{L}}}

\def\cL{{\mathcal{L}}}
\def\cC{{\mathcal{C}}}
\def\cO{\mathcal{O}}
\def\oM{\overline{\mathcal{M}}}
\def\cM{{\mathcal{M}}}
\def\C{{\mathcal{C}}}

\def\Z{\mathbb{Z}}
\def\ZZ{\mathcal{Z}}

\def\C{\mathbb{C}}
\def\Q{\mathbb{Q}}

\def\qed{{\hfill $\Diamond$}}

\def\Aut{{\rm Aut}}

\def\E{\mathrm{E}}
\def\n{\mathrm{n}}
\def\L{\mathrm{L}}
\def\V{\mathrm{V}}
\def\H{\mathrm{H}}
\def\g{\mathrm{g}}
\def\rarr{\rightarrow}

\def\Klog{\omega_{\rm log}}
\def\rmW{\mathsf{W}}
\def\ff{\mathsf{Q}}
\def\fff{\widehat{\mathsf{Q}}}

\def\v{\mathsf{v}}
\def\h{\mathsf{h}}
\def\hH{\mathsf{h}}

\def\bb{\mathsf{m}}
\def\cc{\ell}

\usepackage{theorem}

\newtheorem{theorem}{Theorem}
\newtheorem{prop}{Proposition}

\newtheorem{corollary}[theorem]{Corollary}
\newtheorem{lemma}[prop]{Lemma}

{\theorembodyfont{\rmfamily}

\newtheorem{remark}{Remark}

}

\title{Double ramification cycles on the moduli spaces of curves}
\author{F. Janda, R. Pandharipande, A. Pixton, D. Zvonkine}
\date{March 2016}

\begin{document}

\maketitle

\vspace{-20pt}

\begin{abstract} 
Curves of genus $g$ which admit a map to $\PP^1$
with specified ramification profile $\mu$ over $0\in \PP^1$
and $\nu$ over $\infty \in \PP^1$
define a double ramification cycle $\DR_{g}(\mu,\nu)$ on the moduli 
space of curves.
The study of the restrictions of these cycles to the moduli of
nonsingular curves is a classical
topic. In 2003, Hain calculated the cycles for
curves of compact type. We study here double ramification cycles
on the moduli space of Deligne-Mumford stable curves.

The cycle $\DR_{g}(\mu,\nu)$ for stable curves
is defined via the virtual fundamental class of the moduli of
stable
maps to rubber. Our main result is the 
proof of an explicit formula for $\DR_{g}(\mu,\nu)$
in the tautological ring conjectured by Pixton in 2014. The formula
expresses the double ramification cycle as a sum over stable graphs
(corresponding to strata classes) with summand equal to a
product over markings and edges.
The result answers a question of Eliashberg from 2001
and specializes to Hain's formula
in the compact type case. 

When $\mu=\nu=\emptyset$, the formula for double ramification
cycles expresses
the top Chern class $\lambda_g$ of the Hodge bundle of $\oM_g$
as a push-forward of tautological classes
 supported on the divisor of non-separating nodes. Applications
to Hodge integral calculations are given.

\end{abstract}

\parskip=5pt
\baselineskip=15pt
\pagebreak





\setcounter{section}{-1}
\section{Introduction}

\subsection{Tautological rings}
Let $\oM_{g,n}$ be the moduli space of stable curves
of genus $g$ with $n$ marked points.
Since Mumford's article 
\cite{Mum},
there has been substantial progress in the study of the
structure of the tautological rings{\footnote{All Chow groups
will be taken with $\mathbb{Q}$-coefficients.}} 
$$R^*(\oM_{g,n}) \subset A^*(\oM_{g,n})$$
of the moduli spaces of curves. We refer the reader to \cite{FP-Handbook}
for  
a survey of the basic definitions and properties of $R^*(\oM_{g,n})$.

The recent results \cite{J,PP,PPZ} concerning relations in $R^*(\oM_{g,n})$, conjectured
to be {\em all} relations \cite{P}, may be viewed as parallel
 to the presentation
$$A^*(\mathsf{Gr}(r,n),\Q) = \Q\big[c_1(\mathsf{S}),\ldots,c_r(\mathsf{S})\big]\,
\big/\,\big(s_{n-r+1}(\mathsf{S}), \ldots,
s_n(\mathsf{S})\big)$$
of the Chow ring of the Grassmannian via the Chern and Segre classes of the universal
subbundle
$$S\rightarrow \mathsf{Gr}(r,n)\, .$$
The Schubert calculus for the Grassmannian contains several explicit formulas
for geometric loci. Our main result here is an explicit formula
 for the double ramification cycle in $R^*(\oM_{g,n})$.

\subsection{Double ramification cycles}

\subsubsection{Notation} \label{zz11}
Double ramification data for maps will be specified by a vector
$$A=(a_1,\ldots,a_n)\, , \ \ \ a_i\in \Z$$
 satisfying the balancing condition
$$\sum_{i=1}^n a_i =0\,.$$
The integers $a_i$ are the {\em parts} of $A$.
We separate the positive, negative, and 0 parts of $A$ as follows.
The positive parts of $A$ define a partition
$$\mu = \mu_1+ \ldots+ \mu_{\ell(\mu)}\ .$$
 The negatives of the negative parts of $A$ define a second partition
$$\nu = \nu_1+\ldots+ \nu_{\ell(\nu)}\ .$$
Up to a reordering{\footnote{There will be no
 disadvantage in assuming
$A$ is ordered as in \eqref{fvvf}.}}  of the parts, we have
\begin{equation}\label{fvvf}
A= (\mu_1,\ldots,\mu_{\ell(\mu)}, -\nu_1, \ldots, -\nu_{\ell(\nu)}, 
\underbrace{0,
\ldots, 0}_{n-\ell(\mu)-\ell(\nu)})\, .
\end{equation}

Since the parts of $A$ sum to 0, the partitions $\mu$ and 
$\nu$ must be of the same size. Let
$${{D}}=|\mu| = | \nu|\, $$
be the {\em degree} of $A$.
Let $I$ be the set 
of markings corresponding to
the $0$ parts of $A$.

The various constituents of $A$ are permitted to be empty.
The degree 0 case occurs when
$$\mu=\nu=\emptyset\, .$$
If $I=\emptyset$, then $n=\ell(\mu)+\ell(\nu)$. 
The
empty vector $A$ is permitted, then
   $$n=0\, , D=0\, , \text{ and } \mu=\nu=I=\emptyset\,. $$

\subsubsection{Nonsingular curves}

Let $\cM_{g,n}\subset \oM_{g,n}$ be  
the moduli space of 
nonsingular pointed curves. Let 
$$A=(a_1,\ldots,a_n)$$ be 
a vector of double ramification data as defined in Section \ref{zz11}.
Let $$\ZZ_g(A)\subset \cM_{g,n}$$
be the locus{\footnote{When considering
$\ZZ_g(A)$, we {\em always} assume the stability
condition  $2g-2+n>0$ holds.}} parameterizing curves $[C,p_1,\ldots,p_n]\in \cM_{g,n}$ satisfying
\begin{equation}\label{hh33}
\cO_C\Big(\sum_i{a_i p_i}\Big) \stackrel{\sim}{=} \cO_C\ .
\end{equation}
Condition \eqref{hh33} is algebraic and defines $\ZZ_g(A)$
canonically as a substack of $\cM_{g,n}$.

If $[C,p_1,\ldots,p_n]\in \ZZ_g(A)$, the defining
condition \eqref{hh33}
yields a rational function (up to $\C^*$-scaling),
\begin{equation}\label{vv44}
f:C\rightarrow \PP^1\, ,
\end{equation}
 of degree $D$ with ramification profile $\mu$ over $0\in \PP^1$ and
$\nu$ over $\infty\in \PP^1$.  
The markings corresponding
to $0$ parts lie over $\PP^1 \setminus \{0,\infty\}$.
Conversely, every such morphism \eqref{vv44}, up to $\C^*$-scaling,
determines an element of $\ZZ_g(A)$. 

We may therefore view $\ZZ_g(A)\subset \cM_{g,n}$
as the moduli space of degree $D$ maps (up to $\C^*$-scaling), 
\begin{equation*}
f:C\rightarrow \PP^1\, ,
\end{equation*}
with ramification profiles $\mu$ and $\nu$ over $0$ and $\infty$
respectively. The term {\em double ramification} is motivated
by the geometry of the map $f$.

The dimension of $\ZZ_g(A)$ is easily calculated via
the theory of Hurwitz covers.
Every map \eqref{vv44} can be deformed within $\ZZ_g(A)$
to a map with only {\em simple ramification} over
$\PP^1 \setminus \{0,\infty\}$. The number of simple branch points
in $\PP^1 \setminus \{0,\infty\}$ is determined by
the Riemann-Hurwitz formula to equal
$$\ell(\mu)+\ell(\nu)+2g-2\, .$$
The dimension of the irreducible components of $\ZZ_g(A)$
can be calculated by varying the branch points~\cite{R} to
equal
$$\ell(\mu)+\ell(\nu)+2g-2 + \ell(I)-1 = 2g-3+n\ .$$
The $-1$ on the left is the effect of the $\C^*$-scaling.
Hence, $\ZZ_g(A)$ is of {\em pure codimension} $g$ in $\cM_{g,n}$.


\subsubsection{The Abel-Jacobi map} \label{abj}
Let $\cM_{g,n}^{\ct} \subset \oM_{g,n}$
be the moduli space of curves of compact type. The universal Jacobian
$$
\text{Jac}_{g,n} \to \cM_{g,n}^{\ct}
$$
is the moduli space of line bundles on the universal curve
of  degree 0 on {\em every} irreducible component of {\em every} fiber.

There is a natural compactification of $\ZZ_g(A)$ in the moduli space of compact type curves, 
$$
\ZZ_g(A) \subset \ZZ^{\ct}_g(A)\subset \cM_{g,n}^{\ct}\, ,
$$
obtained from the geometry of the universal Jacobian.

Consider the universal curve over the compact type moduli space
$$\pi^{\ct}: \cC^{\ct} \rightarrow \cM_{g,n}^{\ct}\, .$$
The double ramification vector $A=(a_1,\ldots,a_n)$ determines a line bundle
$\mathsff$ on $\cC^{\ct}$ of relative degree 0,
$$\mathsff = \sum_i a_i [\mathsf{p}_i]\, ,$$
where $\mathsf{p}_i \subset \cC^{\ct}$ is the section of $\pi^{\ct}$ 
corresponding to the marking $p_i$.

By twisting $\mathsff$ by {\em components}
 of the $\pi^{\ct}$ inverse images of the
boundary divisors of $\cM_{g,n}^{\ct}$, we can easily construct a line bundle
${\mathsff}'$ which has degree 0 on {every} component of 
{every} fiber of $\pi^{\ct}$. The result
${\mathsff}'$ is unique up to
twisting by the $\pi^{\ct}$ inverse
images of the boundary divisors of $\cM_{g,n}^{\ct}$.

Via ${\mathsff}'$, we obtain a section of the universal Jacobian,
\begin{equation}\label{ffrr}
\phi: \cM_{g,n}^{\ct} \rightarrow \text{Jac}_{g,n} .
\end{equation}
Certainly the closure of $\ZZ_g(A)$ in $\cM_{g,n}^{\ct}$ lies in the 
$\phi$ inverse image of the $0$-section $S\subset \text{Jac}_{g,n}$
of the relative Jacobian.
We define
\begin{equation}\label{cctt}
\ZZ^{\ct}_g(A) = \phi^{-1}(S) \subset \cM_{g,n}^{\ct}\ .
\end{equation}
The substack $\ZZ^{\ct}_g(A)$ is independent of
the choice of ${\mathsff}'$.
From the definitions, we see
$$\ZZ^{\ct}_g(A)\cap \cM_{g,n}= \ZZ_g(A)\ . $$

While $S\subset \text{Jac}_{g,n}$ is of pure codimension $g$,
$\ZZ_g(A)^{\ct}$ typically has components of excess dimension
(obtained
from genus 0 components of the domain).{\footnote{Hence, 
$\ZZ_g(A)\subset \ZZ^{\ct}_g(A)$ is {\em not} always dense.}}
A cycle class of the expected dimension
$$\DR_g^{\ct}(A)=\left[\ZZ_g^{\ct}(A)\right]^{\vir} \in A^g(\cM_{g,n}^{\ct})$$
is defined by $\phi^{*}([S])$. A closed formula for 
$\DR_g^{\ct}(A)$ was obtained by Hain \cite{Hain}. A simpler approach was later provided
by Grushevsky and Zakharov~\cite{GruZak}.

\subsubsection{Stable maps to rubber} 
Let $A=(a_1,\ldots,a_n)$ be a vector of double ramification data.
The moduli space 
 $$\oM_{g,I}(\PP^1,\mu,\nu)^{\sim}$$ parameterizes 
{\em stable relative maps} of connected curves of genus $g$
to rubber with
ramification profiles $\mu,\nu$ over the point
 $0,\infty \in \PP^1$ respectively.
The tilde indicates a rubber target. We refer the reader to \cite{JLi}
for the foundations.
There is a natural morphism
$$ \epsilon: \oM_{g,I}(\PP^1,\mu,\nu)^{\sim} \rightarrow \oM_{g,n}$$
forgetting everything except the marked domain curve.

Let $R/\{0,\infty\}$ be a rubber target{\footnote{$R$ is a chain of $\PP^1$s.}}, and let
$$[f:C \rightarrow R/\{0,\infty\}]\in {\oM}_{g,I}(\PP^1,\mu,\nu)^{\sim} $$
be a moduli point.
If $C$ is of compact type and $R\stackrel{\sim}{=}\PP^1$, then
$$\epsilon([f])\in \ZZ^{\ct}_g(A)\subset \oM_{g,n}\, . $$
Hence, we have the inclusion
$$  \ZZ^{\ct}_g(A) \subset \text{Im}(\epsilon) \subset \oM_{g,n}\, .$$

The virtual dimension of $\oM_{g,I}(\PP^1,\mu,\nu)^\sim$ is $2g-3+n$. 
We define the {\em double ramification cycle} to be the push-forward
$$\DR_g(A) = \epsilon_*\Big[\oM_{g,I}(\PP^1,\mu,\nu)^\sim\Big]^{\vir} \in A^g(\oM_{g,n})\ .$$

\subsubsection{Basic properties}

The first property of $\DR_g(A)$ is a compatibility
with the Abel-Jacobi construction of Section \ref{abj}.
Let $$\iota: \cM^{\ct}_{g,n} \rightarrow \oM_{g,n}\ $$
be the open inclusion map.{\footnote{Compatibility on the
smaller open set of curves with rational tails was proven
earlier in \cite{Cav}.}}

\begin{prop}[Marcus-Wise \cite{MW}] We have
$$\iota^* \DR_g(A) = \DR_g^{\ct}(A) \in A^g(\cM^{\ct}_{g,n})\ .$$
\end{prop}

By definition, $\DR_g(A)$ is a Chow class on the
moduli space of stable curves.
In fact, the double ramification cycle lies in the tautological ring.

\begin{prop}[Faber-Pandharipande \cite{FP}]  \label{prfp}
$\ \DR_g(A)\in R^g(\oM_{g,n})$.
\end{prop}

The proof of \cite{FP} provides an algorithm to calculate $\DR_g(A)$
in the tautological ring, but the complexity of the method is too great:
there is no apparent way to obtain
an explicit formula for  $\DR_g(A)$ directly from \cite{FP}.

\subsection{Stable graphs and strata}
\subsubsection{Summation over stable graphs}
The strata of $\oM_{g,n}$ are the quasi-projective subvarieties parameterizing pointed curves of
a fixed {\em topological type}. The moduli space $\oM_{g,n}$ is a disjoint union of finitely
many strata. 

The main result of the paper is a proof of an explicit formula conjectured by Pixton~\cite{PixDR} in 2014
for $\DR_g(A)$ in the tautological ring $R^*(\oM_{g,n})$.
The formula is written in terms of a summation over stable graphs $\Gamma$ indexing
the strata of $\oM_{g,n}$. 
We review here the standard indexing of the strata
of $\oM_{g,n}$ by stable graphs.

\subsubsection{Stable graphs}
The strata of the moduli space of curves correspond
to {\em stable graphs} 
$$\Gamma=(\V, \H,\L, \ \mathrm{g}:\V \rarr \Z_{\geq 0},
\ v:\H\rarr \V, 
\ \iota : \H\rarr \H)$$
satisfying the following properties:
\begin{enumerate}
\item[(i)] $\V$ is a vertex set with a genus function $\g:V\to \Z_{\geq 0}$,
\item[(ii)] $\H$ is a half-edge set equipped with a 
vertex assignment $v:\H \to \V$ and an involution $\iota$,
\item[(iii)] $\E$, the edge set, is defined by the
2-cycles of $\iota$ in $\H$ (self-edges at vertices
are permitted),
\item[(iv)] $\L$, the set of legs, is defined by the fixed points of $\iota$ and is
placed in bijective correspondence with a set of markings,
\item[(v)] the pair $(\V,\E)$ defines a {\em connected} graph,
\item[(vi)] for each vertex $v$, the stability condition holds:
$$2\g(v)-2+ \n(v) >0,$$
where $\n(v)$ is the valence of $\Gamma$ at $v$ including 
both half-edges and legs.
\end{enumerate}
An automorphism of $\Gamma$ consists of automorphisms
of the sets $\V$ and $\H$ which leave invariant the
structures $\L$, $\mathrm{g}$, $v$, and $\iota$.
Let $\Aut(\Gamma)$ denote the automorphism group of $\Gamma$.

The genus of a stable graph $\Gamma$ is defined by:
$$\g(\Gamma)= \sum_{v\in V} \g(v) + h^1(\Gamma).$$
A quasi-projective stratum of $\oM_{g,n}$ corresponding to 
 Deligne-Mumford stable curves of fixed topological type
naturally determines
a stable graph of genus $g$ with $n$ legs by considering the dual 
graph of a generic pointed curve parameterized by the stratum.

Let
$\mathsf{G}_{g,n}$
be the set of isomorphism classes of stable graphs of genus $g$ with $n$ legs. 
The set $\mathsf{G}_{g,n}$ is finite.

\subsubsection{Strata classes}
To each stable graph $\Gamma\in \mathsf{G}_{g,n}$, we associate the moduli space
\begin{equation*}
\oM_\Gamma =\prod_{v\in \V} \oM_{\g(v),\n(v)}.
\end{equation*}
 There is a
canonical
morphism 
\begin{equation}\label{dwwd}
\xi_{\Gamma}: \oM_{\Gamma} \rarr \oM_{g,n}
\end{equation}
 with image{\footnote{
The degree of $\xi_\Gamma$ is $|\Aut(\Gamma)|$.}}
equal to the closure of the stratum
associated to the graph $\Gamma$.  To construct $\xi_\Gamma$, 
a family of stable pointed curves over $\oM_\Gamma$ is required.  Such a family
is easily defined 
by attaching the pull-backs of the universal families over each of the 
$\oM_{\g(v),\n(v)}$  along the sections corresponding to the two halves
of each edge in $\E(\Gamma)$.

Each half-edge $h\in \H(\Gamma)$ determines a cotangent line
$$\mathcal{L}_h \rightarrow \oM_\Gamma\, .$$
If $h\in \L(\Gamma)$, then $\mathcal{L}_h$ is the pull-back
via $\xi_\Gamma$ of the corresponding cotangent line of $\oM_{g,n}$.
If $h$ is a side of an edge $e\in \E(\Gamma)$, then 
$\mathcal{L}_h$ is the cotangent line of  the corresponding
side of a node.
Let 
$$\phi_h=c_1(\mathcal{L}_h)\, \in A^1(\oM_\Gamma,\Q)\, .$$

\subsection{Pixton's formula} \label{pixfor}

\subsubsection{Weightings} \label{Sssec:weightings}
Let $A=(a_1,\ldots,a_n)$ be a vector of  double ramification data. Let 
$\Gamma \in \mathsf{G}_{g,n}$ be a stable graph of genus $g$ with $n$ legs.
A {\em weighting}
of $\Gamma$ is a function on the set of half-edges,
$$ w:\H(\Gamma) \rightarrow \Z,$$
which satisfies the following three properties:
\begin{enumerate}
\item[(i)] $\forall h_i\in \L(\Gamma)$, corresponding to
 the marking $i\in \{1,\ldots, n\}$,
$$w(h_i)=a_i\ ,$$
\item[(ii)] $\forall e \in \E(\Gamma)$, corresponding to two half-edges
$h,h' \in \H(\Gamma)$,
$$w(h)+w(h')=0\, ,$$
\item[(iii)] $\forall v\in \V(\Gamma)$,
$$\sum_{v(h)= v} w(h)=0\, ,$$ 
where the sum is taken over {\em all} $\mathsf{n}(v)$ half-edges incident to $v$.
\end{enumerate}
If the graph $\Gamma$ has loops, $\Gamma$ may carry infinitely many admissible weightings.

Let $r$ be a positive integer.
A {\em weighting mod $r$}
of $\Gamma$ is a function,
$$ w:\H(\Gamma) \rightarrow \{0,\ldots, r-1\},$$
which satisfies exactly properties (i-iii) above, but
with the equalities replaced, in each case, by the condition of 
{\em congruence $mod$ $r$}.
For example, for (i), we require
$$w(h_i)=a_i \mod r\, .$$
Let $\mathsf{W}_{\Gamma,r}$ be the set of weightings mod $r$
of $\Gamma$. The set $\mathsf{W}_{\Gamma,r}$ is finite, with cardinality $r^{h^1(\Gamma)}$.
We view $r$ as a {\em regularization parameter}.

\subsubsection{Pixton's conjecture}
\label{pixconj}

Let $A = (a_1, \dots, a_n)$ be a vector of double ramification data. Let $r$ be a positive
integer.
We denote by
$\P_g^{d,r}(A)\in R^d(\oM_{g,n})$ the degree $d$ component of the tautological class 
\begin{multline*}
\hspace{-10pt}\sum_{\Gamma\in \mathsf{G}_{g,n}} 
\sum_{w\in \mathsf{W}_{\Gamma,r}}
\frac1{|\Aut(\Gamma)| }
\, 
\frac1{r^{h^1(\Gamma)}}
\;
\xi_{\Gamma*}\Bigg[
\prod_{i=1}^n \exp(a_i^2 \psi_{h_i}) \cdot 
\\ \hspace{+10pt}
\prod_{e=(h,h')\in \E(\Gamma)}
\frac{1-\exp(-w(h)w(h')(\psi_h+\psi_{h'}))}{\psi_h + \psi_{h'}} \Bigg]\, .
\end{multline*} 
in $R^*(\oM_{g,n})$.

Inside the push-forward in the above formula, the first product 
$$\prod_{i=1}^n \exp(a_i^2 \psi_{h_i})\, $$
is over $h\in L(\Gamma)$
via the correspondence of legs and markings.
The second product is over all $e\in E(\Gamma)$.
The factor 
$$\frac{1-\exp(-w(h)w(h')(\psi_h+\psi_{h'}))}{\psi_h + \psi_{h'}}$$
is well-defined since 
\begin{enumerate}
\item[$\bullet$]
the denominator formally divides
the numerator,
\item[$\bullet$] the factor is symmetric in $h$ and $h'$.
\end{enumerate}
No edge orientation is necessary.

The following fundamental polynomiality property of $\P_g^{d,r}(A)$  has been proven by Pixton.

\begin{prop}[Pixton \cite{PixDR2}] For fixed $g$, $A$, and $d$, the \label{pply}
class
$$\P_g^{d,r}(A) \in R^d(\oM_{g,n})$$
is polynomial in $r$ (for all sufficiently large $r$).
\end{prop}

We denote by $\P_g^d(A)$ the value at $r=0$ 
of the polynomial associated to $\P_g^{d,r}(A)$ by Proposition~\ref{pply}. In other words, $\P_g^d(A)$ is the {\em constant} term of the associated polynomial in $r$. For the reader's convenience, 
Pixton's proof of Proposition \ref{pply} is given in the Appendix.

The main result of the paper is a proof of the 
 formula for double ramification
cycles conjectured earlier by Pixton~\cite{PixDR}.

\begin{theorem} \label{FFFF}
For $g\geq 0$ and double ramification data $A$, we have
$$\DR_g(A) = 2^{-g}\, \P_g^g(A)\, \in R^g(\oM_{g,n}).$$ 
\end{theorem}

For $d<g$,  the classes $\P_g^d(A)$ do not yet have a
geometric interpretation. However
for $d>g$, the following vanishing conjectured by Pixton is now established.

\begin{theorem}[Clader-Janda \cite{cj}] \label{Thm:van}
For all $g\geq 0$, double ramification data $A$, and $d>g$, we have
$$\P_g^d(A) = 0 \, \in R^d(\oM_{g,n}) .$$ 
\end{theorem}

Pixton \cite{PixDR} has further proposed a twist of the
formula for $\P_g^d(A)$ by $k\in \Z$. 
The codimension $g$ class in the
$k=1$ case has been (conjecturally) related to the moduli
spaces of meromorphic differential in the Appendix of \cite{FarP}.
The  vanishing result of Clader-Janda~\cite{cj} proves
Pixton's vanishing conjecture in codimensions $d>g$ for all $k$.
The $k$-twisted theory will be discussed in Section~\ref{Sec:rspin}.

In a forthcoming paper \cite{forth}, we will generalize Pixton's formula to the situation of maps to a ${\mathbb{P}}^1$-rubber bundle over a target manifold~$X$.

\subsection{Basic examples}

\subsubsection{Genus 0} 
Let $g=0$ and 
$A=(a_1,\ldots,a_n)$ be a vector of double ramification data.
The graphs $\Gamma\in\mathsf{G}_{0,n}$ are trees. Each $\Gamma\in \mathsf{G}_{0,n}$ has a {\em unique} weighting mod $r$ (for each $r$). 
We can directly calculate:
\begin{equation} \label{jj12}
\P_0(A) = \sum_{d\geq 0} \P_0^d(A)= \exp \left[ \sum_{i=1}^n a_i^2 \psi_i 
- \frac12 
\sum_{
\substack{I \sqcup J = \{1, \dots, n \} \\
|I|, |J| \geq 2}
}
a_I^2 \delta_{I,J}
\right],
\end{equation}
where $a_I = \sum_{i \in I} a_i$.

By Theorem \ref{FFFF} , the double ramification cycle is the degree~0 term
of \eqref{jj12}, 
$$\DR_0(A) = 2^{-0}\, \P_0^0(A)= 1\, .$$
 The vanishing of Theorem \ref{Thm:van} implies  
$$
\sum_{i=1}^n a_i^2 \psi_i - \frac12 
\sum_{
\substack{I \sqcup J = \{1, \dots, n \} \\
|I|, |J| \geq 2}
}
a_I^2 \delta_{I,J} = 0 \ \in R^1(\oM_{0,n})\ , $$
where $\delta_{I,J}\in R^1(\oM_{0,n})$ is the class of the boundary
divisor with marking distribution $I \sqcup J$.
The above vanishing can be proven here directly by expressing 
$\psi$ classes in terms of boundary divisors.

\subsubsection{Genus 1}
Let $g=1$ and $A=(a_1,\ldots,a_n)$ be a vector of double ramification data.
Denote by 
$\delta_0\in R^1(\oM_{1,n})$ 
the boundary divisor of singular stable curves with a nonseparating node. Our convention is  
$$\delta_0 = \frac12\, \xi_* \Big[\oM_{0,n+2}\Big]\, , \ \ \ \ \  \xi:\oM_{0,n+2} \rightarrow \oM_{1,n}\, .$$ 
The division by the order of the automorphism group
 is included in the definition of $\delta_0$.
Denote by $$\delta_I\in R^1(\oM_{1,n})$$ 
the class of the boundary divisor of curves with a rational component 
carrying the markings $I \subset \{1, \dots, n\}$ and an 
elliptic component carrying the markings $I^c$.  

We can calculate directly from the definitions:
\begin{equation}\label{yy99}
\P_1^1(A) = \sum_{i=1}^n a_i^2 \psi_i - 
\sum_{
\substack{I \subset \{1, \dots, n \} \\
|I| \geq 2}
}
a_I^2 \delta_I - \frac16 \delta_0\, .
\end{equation}
As before, $a_I = \sum_{i \in I} a_i$.
By Theorem \ref{FFFF}, 
$$\DR_1(A)=\frac{1}{2}\, \P_1^1(A)\, .$$ 
In genus $1$, the virtual class plays essentially no role
(since the moduli spaces of maps to rubber are of
expected dimension).
The genus 1 formula was already known to Hain \cite{Hain}.

It is instructive to compute the coefficient $-\frac{1}{6}$ of $\delta_0$. 
The class $\delta_0$ corresponds to the graph $\Gamma$ with one vertex of 
genus~0 and one loop,  $$h^1(\Gamma) = 1\, .$$  
According to the definition, the coefficient of $\delta_0$ is the $r$-free term of the polynomial
$$
\frac1r \left[ \sum_{w=1}^{r-1} w(r-w) \right],
$$
which is $-B_2 = -\frac{1}{6}$.

\subsubsection{Degree 0} 

Let $g\geq 0$ and $A=(0,\ldots,0)$, so $\mu=\nu=\emptyset$.
The canonical map $\epsilon$ from  
the moduli of stable maps to the moduli of curves is then
an isomorphism
$$ \epsilon: \oM_{g,n}(\PP^1,\emptyset,\emptyset)^{\sim} 
\stackrel{\sim}{\longrightarrow}
 \oM_{g,n}\ .$$
By an analysis of the obstruction theory,
$$\DR_g(0,\ldots,0) = (-1)^g\lambda_g \ \in R^g(\oM_{g,n})\ $$
 where $\lambda_g$ is the top Chern class of the Hodge bundle
$$\mathbb{E} \rightarrow \oM_{g,n}\ . $$

Let $\Gamma\in \mathsf{G}_{g,n}$ be a stable graph. By the definition
of a {\em weighting mod $r$} and the condition $a_i=0$ for all markings $i$,
the weights 
$$w(h)\, ,\ w(h')$$  on the two halves of {\em every} separating 
edge $e$ of $\Gamma$ must both be $0$. The factor in Pixton's formula for 
$e$,
$$\frac{1-\exp(-w(h)w(h')(\psi_h+\psi_{h'}))}{\psi_h + \psi_{h'}}\, ,$$
then vanishes and kills the contribution of $\Gamma$ to $\P_g^g(0,\ldots,0)$.

Let $g\geq 1$.
Let $\xi$ be the map to $\oM_{g,n}$ associated to the divisor
of curves with non-separating nodes
$$\xi:\oM_{g-1,n+2} \rightarrow \oM_{g,n}\, .$$
Since graphs $\Gamma$ with separating edges do not contribute to 
$\P_g^g(0,\ldots,0)$
and the trivial graph with no edges does not contribute in codimension $g$,
$$ \P_g^g(0,\ldots,0) = \xi_*\, \Lambda_{g-1}^{g-1}(0,\ldots,0)$$
for an explicit tautological class
$$\Lambda_{g-1}^{g-1}(0,\ldots,0) \ \in R^{g-1}(\oM_{g-1,n+2})\ .$$

\begin{corollary} \label{xzz} For $g\geq 1$, we have
$$\lambda_g = (-2)^{-g}\, \xi_*\, \Lambda_{g-1}^{g-1}(0,\ldots,0)\ 
\in R^g(\oM_{g,n})\ .$$
\end{corollary}

The class $\lambda_g$ is easily proven to vanish on $\cM_{g,n}^{\ct}$
via the Abel-Jacobi map \eqref{ffrr} and well-known properties of the moduli
space of abelian varieties \cite{Geer}. Hence, there exists a {\em Chow} class
$\gamma_{g-1,n+2} \in 
A^{g-1}(\oM_{g-1,n+2})$ satisfying
$$\lambda_g = \xi_* \gamma_{g-1,n+2}\ \in A^g(\oM_{g,n})\, .$$
Corollary \ref{xzz} is a much stronger statement since
$\Lambda_{g-1}^{g-1}(0,\ldots,0)$ is a {\em tautological} class
on $\oM_{g-1,n+2}$ given by an explicit formula.
No such expressions were known before.

\vspace{12pt}
\noindent $\bullet$ For $\oM_{1,n}$, Corollary \ref{xzz} yields
$$\lambda_1 = -\frac{1}{2}\, \xi_*\, \Lambda_{1-1}^{1-1}(0,\ldots,0)\ \in R^1(\oM_{1,n})\, .$$
Moreover, by  calculation \eqref{yy99},
$$\Lambda_{1-1}^{1-1}(0,\ldots, 0) = -\frac{1}{12}\, \Big[\oM_{0,n+2}\Big]\, ,$$
since $\delta_0$ already carries a $\frac{1}{2}$ for the graph automorphism.
Hence, we recover the well-known formula
$$\lambda_1 = \frac{1}{24}\, \xi_*\Big[\oM_{0,n+2}\Big]\ \in R^1(\oM_{1,n}) \, .$$

\vspace{12pt}
\noindent $\bullet$ For $\oM_{2}$, Corollary \ref{xzz} yields
$$\lambda_2 = \frac{1}{4}\, \xi_* \Lambda_{2-1}^{2-1}(\emptyset)\ 
\in R^2(\oM_{2})\, .$$
Denote by $\alpha$ the class of the boundary stratum of genus~0 curves with two self-intersections.
We include the division by the order of the automorphism group in the
 definition, 
$$\alpha = \frac18\ \xi_* \Big[\oM_{0,4}\Big]\ .$$ 
Denote by $\beta$ the class of the boundary divisor of genus~1 curves with one self-intersection multiplied 
by the $\psi$-class at one of the branches. Then we have
\begin{equation}\label{wqq}
\xi_*\Lambda_{2-1}^{2-1}(\emptyset) = \frac1{36} \alpha + \frac1{60} \beta\, .
\end{equation}
The classes $\alpha$ and $\beta$ form a basis of $R^2(\oM_2)$. The expression 
 $$\lambda_2 = \frac{1}{144} \alpha + \frac{1}{240} \beta $$
 can be checked by intersection with the two boundary divisors of $\oM_2$.

The coefficient of $\alpha$ in \eqref{wqq} is obtained from the $r$-free term of
$$
\frac1{r^2} \sum_{1 \leq w_1, w_2 \leq r-1} w_1 w_2 (r-w_1) (r-w_2)\, .
$$
The answer is $B_2^2=\frac{1}{36}$.
The coefficient of $\beta$ is the $r$-free term of
$$
-  \frac{1}{2r} \sum_{w=1}^{r-1} w^2(r-w)^2
$$
given by $-\frac{1}{2}B_4 = \frac{1}{60}$.


\subsection{Strategy of proof}

The proof of Proposition \ref{prfp} given in \cite{FP} uses the Gromov-Witten theory of
$\PP^1$ relative to the point $\infty\in \PP^1$. The localization
relations of \cite{FP} intertwine Hodge classes
on the moduli spaces of curves with cycles obtained from maps to rubber. The
interplay is very complicated with many auxiliary classes (not just the double ramification
cycles).  

Our proof of Theorem \ref{FFFF} is obtained by studying the Gromov-Witten
theory of $\PP^1$ with an orbifold{\footnote{$B\Z_r$ is
the quotient stack obtained by the trivial action of 
$\Z_r\cong\frac{\Z}{r\Z}$ on a point.}} point $B\Z_r$ at $0\in \PP^1$ 
and a relative 
point  $\infty\in \PP^1$. Let
$(\PP^1[r],\infty)$ denote the resulting orbifold/relative geometry.{\footnote{The target $(\PP^1[r],\infty)$ was previously studied
in \cite{JPT} in the context of Hurwitz-Hodge integrals.}}  
The proof may be explained conceptually as studying the localization relations 
for all $(\PP^1[r],\infty)$  {\em simultaneously}. For large $r$, the localization relations
(after appropriate rescaling)
may be viewed as having polynomial coefficients in $r$. Then the $r=0$ specialization 
of the linear algebra exactly yields Theorem \ref{FFFF}.  The regularization parameter $r$ in the definition of
Pixton's formula in Section \ref{pixfor} is related to the orbifold structure $B\Z_r$.

Our proof is inspired by the
structure of Pixton's formula --- especially  the existence of the regularization parameter $r$.
The argument requires two main steps:

\noindent $\bullet$ Let $\oM^r_{g;a_1, \dots, a_n}$ be the moduli space of tensor 
$r$th roots of 
the bundles $\cO_C(\sum_{i=1}^n a_i x_i)$ associated to pointed curves $[C,p_1,\ldots,p_n]$. 
Let 
$$\pi : \cC^r_{g;a_1, \dots, a_n} \to \oM^r_{g;a_1, \dots, a_n}$$ 
be the universal curve, and let $\cL$ be the universal $r$th root over the universal curve. 
Chiodo's formulas \cite{Chiodo2} allow us to prove 
that the push-forward of $r\cdot c_g(-R^*\pi_*\cL)$ to $\oM_{g,n}$ and $2^{-g}\,
\P_g^{g,r}(A)$  
are polynomials in 
$r$ with the {\em same} constant term. 
Pixton's formula therefore has a geometric interpretation in terms of the intersection theory of 
the moduli space of $r$th roots.

\noindent $\bullet$
 We use the localization formula \cite{GrP} for the virtual
class of the moduli space of stable maps to the orbifold/relative geometry 
$(\PP^1[r],\infty)$.
The monodromy conditions at $0$ are given by $\mu$ and the ramification profile over $\infty$ is given by 
$\nu$.
The push-forward of the localization formula to $\oM_{g,n}$ is a Laurent series in the equivariant parameter $t$
and in $r$. The coefficient of $t^{-1}r^0$  must vanish by geometric considerations. 
We prove the relation obtained from
the coefficient of $t^{-1}r^0$  has {\em only two terms}.
  The first is the
 constant term in $r$ of the push-forward of 
$r\cdot c_g(-R^*\pi_*\cL)$ to $\oM_{g,n}$. 
The second term is the double ramification cycle 
$\DR_g(A)$
with a minus sign. The vanishing of
the sum of the two terms yields Theorem \ref{FFFF}.

\subsection{Acknowledgments} 
We are grateful to A.~Buryak,
R.~Cavalieri, A.~Chiodo, E.~Clader, C.~Faber, G.~Farkas,
S.~Grushevsky, J.~Gu\'er\'e, G.~Oberdieck, D.~Oprea, R.~Vakil, and Q.~Yin for 
 discussions about
double ramification cycles and related topics.
The {\em Workshop on Pixton's conjectures} at
ETH Z\"urich  in October 2014 played in an important
role in our collaboration. A.~P. and D.~Z. have been
frequent guests of the 
Forschungsinstitut f\"ur Mathematik (FIM) at ETHZ.
The paper was completed at the conference {\em Moduli spaces of
holomorphic differentials} at Humboldt University in  Berlin 
in February 2016 attended
by all four authors.

F.~J. was partially supported by the Swiss National Science Foundation
grant SNF-200021143274. 
R.~P. was partially supported by 
 SNF-200021143\-274, SNF-200020162928, ERC-2012-AdG-320368-MCSK, SwissMap, and
the Einstein Stiftung. 
A.~P. was supported by a fellowship from the Clay Mathematics Institute.
D.~Z. was supported by the grant ANR-09-JCJC-0104-01.

\section{Moduli of stable maps to $B\Z_r$}
\label{Sec:rspin}

\subsection{Twisting by $k$}
Let $k\in \Z$. A vector $A=(a_1,\ldots, a_n)$
of {\em $k$-twisted double ramification data} for genus $g\geq 0$
is defined by the condition{\footnote{We  {\em always} assume the stability
condition  $2g-2+n>0$ holds.}}
\begin{equation}\label{tqqt}
\sum_{i=1}^n a_i \ = \ k(2g-2+n)\ .
\end{equation}
For $k=0$, the condition \eqref{tqqt} does {\em not} depend upon
on the genus $g$ and specializes to the balancing condition of
the double ramification data of Section \ref{zz11}.

We recall the definition \cite{PixDR} of Pixton's $k$-twisted cycle 
$\P_{g}^{d,k}(A)\in R^d(\oM_{g,n})$. In case $k=0$,
$$\P_{g}^{d,0}(A)= \P_{g}^d(A) \in R^d(\oM_{g,n})\ .$$

Let 
$\Gamma \in \mathsf{G}_{g,n}$ be a stable graph of genus $g$ with $n$ legs.
A {\em $k$-weighting mod r}
of $\Gamma$ is a function on the set of half-edges,
$$ w:\H(\Gamma) \rightarrow \{ 0,\ldots, r-1\},$$
which satisfies the following three properties:
\begin{enumerate}
\item[(i)] $\forall h_i\in \L(\Gamma)$, corresponding to
 the marking $i\in \{1,\ldots, n\}$,
$$w(h_i)=a_i \mod r \, ,$$
\item[(ii)] $\forall e \in \E(\Gamma)$, corresponding to two half-edges
$h,h' \in \H(\Gamma)$,
$$w(h)+w(h')=0 \mod r\,$$
\item[(iii)] $\forall v\in \V(\Gamma)$,
$$\sum_{v(h)= v} w(h)=k\big( 2\g(v)-2+\n(v)\big) \mod r\, ,$$ 
where the sum is taken over {\em all} $\mathsf{n}(v)$ half-edges incident to $v$.
\end{enumerate}
We denote by $\rmW_{\Gamma,r,k}$ the set of all $k$-weightings mod $r$ of~$\Gamma$.

Let $A$ be a vector of $k$-twisted ramification data for genus $g$. 
For each positive integer $r$,  let
$\P_g^{d,r,k}(A)\in R^d(\oM_{g,n})$ be the {\em degree $d$} component
 of the tautological class 
\begin{multline*}
\hspace{-10pt}\sum_{\Gamma\in \mathsf{G}_{g,n}} 
\sum_{w\in \mathsf{W}_{\Gamma,r,k}}
\frac1{|\Aut(\Gamma)| }
\, 
\frac1{r^{h^1(\Gamma)}}
\;
\xi_{\Gamma*}\Bigg[ \prod_{v \in \V(\Gamma)} e^{-k^2 \kappa_1(v)}
\prod_{i=1}^n e^{a_i^2 \psi_{h_i}} \cdot 
\\ \hspace{+10pt}
\prod_{e=(h,h')\in \E(\Gamma)}
\frac{1-e^{-w(h)w(h')(\psi_h+\psi_{h'})}}{\psi_h + \psi_{h'}} \Bigg]\, .
\end{multline*} 
in $R^*(\oM_{g,n})$.

All the conventions here are the same as in the definitions of Section 
\ref{pixconj}. 
A new factor appears: $\kappa_1(v)$ is the first $\kappa$ class{\footnote{Our
convention is 
$\kappa_i= \pi_*(\psi_{\n(v)+1}^{i+1})\ \in R^{i}(\oM_{\g(v),\n(v)})$ where
$$\pi: \oM_{\g(v),\n(v)+1} \rightarrow \oM_{\g(v),\n(v)}$$ is the
map forgetting the marking $\n(v)+1$. For a review of $\kappa$ and
and cotangent $\psi$ classes, see \cite{GrP2}.}},
$$\kappa_1(v) \in R^1(\oM_{\g(v),\n(v)})\, ,$$ 
on the moduli space corresponding to a vertex $v$ of the graph $\Gamma$.
As in the untwisted case, polynomiality holds (see the Appendix).

\vspace{8pt}
\noindent {\bf{Proposition $\bf{3'}$} (Pixton \cite{PixDR2})} {\em For fixed $g$, $A$, $k$, and $d$, the 
class
$$\P_g^{d,r,k}(A) \in R^d(\oM_{g,n})$$
is polynomial in $r$ (for all sufficiently large $r$).}

\vspace{8pt}
We denote by $\P_g^{d,k}(A)$ the value at $r=0$ 
of the polynomial associated to $\P_g^{d,r,k}(A)$ by Proposition $3'$. In other words, $\P_g^{d,k}(A)$ is the {\em constant} term of the associated polynomial in $r$. 

In case $k=0$, Proposition $3'$ specializes to Proposition \ref{pply}. Pixton's proof
is given in the Appendix.

\subsection{Generalized $r$-spin structures}
\label{Ssec:koverr}

\subsubsection{Overview}
Let $[C, p_1, \dots, p_n]\in \cM_{g,n}$ be a nonsingular
 curve with distinct  markings and canonical bundle $\omega_C$.
Let
 $$\Klog = \omega_C(p_1 + \dots + p_n)$$ be
the {\em log canonical} bundle.

Let $k$ and $a_1, \dots, a_n$ be integers for which $k(2g-2+n)-\sum a_i$ is divisible by a
positive integer
$r$. 
Then $r$th roots $L$ of the line bundle $$\Klog^{\otimes k}\left(-\sum 
a_i p_i\right)$$ on $C$ 
{\em exist}. 
The space of such $r$th tensor roots possesses a natural compactification 
$\oM_{g;a_1, \dots, a_n}^{r,k}$ constructed in~\cite{Chiodo,Jarvis}. 
The moduli space of $r$th roots carries a universal curve 
$$\pi : \cC_{g;a_1, \dots, a_n}^{r,k} \to \oM_{g;a_1, \dots, a_n}^{r,k}$$ and a universal line bundle $$\cL \to \cC_{g;a_1, \dots, a_n}^{r,k}$$ equipped with an isomorphism
\begin{equation}\label{dd22}
\cL^{\otimes r} \stackrel{\sim}\longrightarrow \Klog^{\otimes k}\left(-\sum 
a_i x_i\right)
\end{equation}
 on $\cC^{r,k}_{g;a_1,\ldots,a_n}$.
In case $k=0$, $\oM_{g;a_1, \dots, a_n}^{r,0}$ is the space of stable maps \cite{AGV}
to the orbifold $B\Z_r$. 

We review here the basic properties of the moduli spaces
$\oM_{g;a_1, \dots, a_n}^{r,k}$ which we require and refer the reader
to \cite{Chiodo, Jarvis} for a foundational development of the theory.

\subsubsection{Covering the moduli of curves.} \label{cmc}
The forgetful map to the moduli of curves 
$$\epsilon:\oM_{g;a_1, \dots, a_n}^{r,k} \to \oM_{g,n}$$
 is a $r^{2g}$-sheeted covering ramified over the boundary divisors. The degree of $\epsilon$ equals $r^{2g-1}$ because each $r$th root 
$L$ has a $\Z_r$ symmetry group 
obtained from multiplication by 
$r$th roots of unity in the fibers.{\footnote{The $\Z_r$-action
respects the 
isomorphism \eqref{dd22}.}}

To transform $\epsilon$ to an {\em unramified} covering in the orbifold sense, the orbifold structure of $\oM_{g,n}$ must be altered. 
We introduce an extra $\Z_r$ stabilizer for each node of each stable curve, see~\cite{Chiodo}. The new orbifold thus obtained is called the 
moduli space of {\em $r$-stable curves}.

\subsubsection{Boundary strata.}
The $r$th roots of the trivial line bundle over the cylinder are classified by a pair of integers 
\begin{equation}\label{cxx3}
0\leq a',a''\leq r-1\, , \ \ \ \ a' + a'' = 0 \mod r
\end{equation}
attached to the borders of the cylinder. 
Therefore, the nodes of stable curves in the compactification 
 $\oM_{g;a_1, \dots, a_n}^{r,k}$ carry the data \eqref{cxx3} with 
 $a'$ and $a''$ assigned to the two branches meeting at the node. 

In case  the node is {\em separating}, the data~\eqref{cxx3}
 are determined uniquely by the way the genus and the marked points 
are distributed over the two components of the curve. In  case 
the node is {\em nonseparating}, 
the data can be arbitrary and describes different boundary divisors. 
 
More generally, the boundary strata of $\oM_{g;a_1, \dots,a_n}^{r,k}$ are described by stable graphs $\Gamma\in \mathsf{G}_{g,n}$ together
a $k$-weighting mod~$r$. The data \eqref{cxx3} is obtained
from the weighting $w$ of the half-edges of $\Gamma$.

\subsubsection{Normal bundles of boundary divisors.}
\label{Sssec:normalbundle}
The normal bundle{\footnote{More precisely, the normal bundle of the associated gluing morphism $\xi$.}
of a boundary divisor in $\oM_{g,n}$ has first 
Chern class $$-(\psi'+ \psi'')\ ,$$ where $\psi'$ and $\psi''$ are the $\psi$-classes associated to the branches of the node. 
In the moduli space of $r$-stable curves {\em and} 
in the space $\oM_{g;a_1, \dots, a_n}^{r,k}$, 
the normal bundle to a boundary divisor has first Chern class 
\begin{equation} \label{Eq:normal}
-\frac{\psi'+ \psi''}{r}\ ,
\end{equation} 
where the classes $\psi'$ and $\psi''$ are pulled back from $\oM_{g,n}$. 

The two standard boundary maps 
$$
\oM_{g-1,n+2} \stackrel{i} \to \oM_{g,n}
$$
and
$$
\oM_{g_1, n_1+1} \times \oM_{g_2, n_2+1} \stackrel{j} \to \oM_{g,n}
$$
are, in the space of $r$th roots, replaced by diagrams
\begin{equation} \label{Eq:phi}
 \oM^{r,k}_{g-1; a_1, \dots, a_n, a', a''} \stackrel{\phi}\longleftarrow \mathcal{B}^{\rm nonsep}_{a',a''} \stackrel{i}\longrightarrow \oM^{r,k}_{g; a_1, \dots, a_n}
\end{equation}
and
\begin{equation} \label{Eq:varphi}
\oM^{r,k}_{g_1; a_1, \dots, a_{n_1}, a'} \times  \oM^{r,k}_{g_2; a'',a_{n_1+1}, \dots, a_n} \stackrel{\varphi}\longleftarrow \mathcal{B}^{\rm sep}_{a',a''} \stackrel{j}\longrightarrow  \oM^{r,k}_{g; a_1, \dots, a_n}.
\end{equation}
The maps $\phi$ and $\varphi$ are of degree~1, but in general are not isomorphisms, see~\cite{ChiZvo}, Section~2.3.

\begin{remark} \label{Rem:2conventionsA}
There are two natural ways to introduce the orbifold structure on the space of $r$th roots. In the first, every node of the curve contributes an extra $\Z/r\Z$ in the stabilizer. In the second, a node of type $(a',a'')$ contributes an extra $\Z/{\rm gcd}(a',r) \Z$ to the stabilizer. The first is used in~\cite{Chiodo}, the second in~\cite{AGV}. 
Here, 
we follow the first convention to avoid gcd factors appearing the computations. 
In Remark~\ref{Rem:2conventionsB}, we indicate the changes which must be
made in the computations if the second convention is follow.  
The final result is the {\em same} (independent of the choice of convention).
\end{remark}

\subsubsection{Intersecting boundary strata.} Let $(\Gamma_A, w_A)$ and $(\Gamma_B, w_B)$ be two stable graphs in $\mathsf{G}_{g,n}$
equipped with $k$-weightings mod~$r$. The intersection of the two corresponding boundary strata of $\oM_{g;a_1, \dots,a_n}^{r,k}$
 can be determined as follows.{\footnote {See the Appendix of \cite{GrP2}
for a  detailed discussion in the case of 
$\oM_{g,n}$.}} Enumerate all pairs $(\Gamma,w)$ of
stable graphs with $k$-weightings mod~$r$ whose edges are marked with letters $A$, $B$ or both in such a way that contacting all edges outside $A$ yields
 $(\Gamma_A, w_A)$ and contracting all edges outside $B$ yields $(\Gamma_B, w_B)$. Assign a factor 
$$-\frac{1}{r}(\psi'+ \psi'')$$ to every edge that is marked by both $A$ and $B$. The sum of boundary strata given by these graphs multiplied by the appropriate combinations of $\psi$-classes represents the intersection of the boundary components $(\Gamma_A, w_A)$ and $(\Gamma_B, w_B)$.

\subsection{Chiodo's formula} \label{Ssec:Chiodo}
Let $g$ be the genus, and let $k\in \Z$ and 
$a_1, \dots, a_n\in \Z$ satisfy
 $$k(2g-2+n)-\sum a_i = 0  \mod r$$ 
for an integer $r>0$. 

Following the notation of Section \ref{Ssec:koverr},  let
$$\pi : \cC_{g;a_1, \dots, a_n}^{r,k} \to \oM_{g;a_1, \dots, a_n}^{r,k}$$
 be the universal curve and $\cL \to \cC_{g;a_1, \dots, a_n}^{r,k}$ the universal 
$r$th root. Our aim here is to describe the total Chern class 
$$c\big(-R^*\pi_*\cL\big)\in A^*(\oM_{g;a_1, \dots, a_n}^{r,k},\Q)$$ 
following the work of A.~Chiodo \cite{Chiodo2}.

Denote by $B_m(x)$ the $m$-th Bernoulli polynomial,
$$\sum_{n=0}^\infty B_n(x)\frac{t^n}{n!}= \frac{te^{xt}}{e^t-1}\, .$$

\begin{prop} \label{Prop:ChiodoExp}
The total Chern class
$c\big(-R^*\pi_*\cL\big)$
is equal to 
\begin{multline*}
\hspace{-10pt}\sum_{\Gamma\in \mathsf{G}_{g,n}} 
\sum_{w\in \mathsf{W}_{\Gamma,r,k}}
\frac1{|\Aut(\Gamma)| }
\, 
r^{|\E(\Gamma)|}
\;
(\xi_{\Gamma,w})_*\Bigg[ \prod_{v \in \V(\Gamma)} e^{- \sum\limits_{m\geq 1} (-1)^{m-1} \frac{B_{m+1}(k/r)}{m(m+1)}\kappa_m(v)} \, \cdot\\
\hspace{+10pt}
\prod_{i=1}^n e^{\sum\limits_{m\geq 1}(-1)^{m-1}\frac{B_{m+1}(a_i/r)}{m(m+1)} \psi^m_{h_i}}\; \cdot 
\prod_{
\substack{e\in \E(\Gamma) \\ e = (h,h')}
}   
\frac{1-e^{\sum\limits_{m \geq 1} (-1)^{m-1}\frac{B_{m+1}(w(h)/r)}{m(m+1)} [(\psi_h)^m-(-\psi_{h'})^m]}}{\psi_h + \psi_{h'}} \Bigg]\, 
\end{multline*} 
in $A^*(\oM_{g;a_1, \dots, a_n}^{r,k},\Q)$.
\end{prop}

Here, $\xi_{\Gamma,w}$ is the map of the boundary stratum corresponding to 
$(\Gamma,w)$ to the moduli space $\oM_{g;a_1, \dots, a_n}^{r,k}$,
$$\xi_{\Gamma,w}: \oM_{\Gamma,w} \rightarrow \oM_{g;a_1,\ldots,a_n}^{r,k}\, .$$
 Every factor in the last product is symmetric in $h$ and $h'$ since
 $$B_{m+1}\left(\frac{w}{r}\right) = (-1)^{m+1}B_{m+1}\left(\frac{r-w}{r}\right)\, $$
for $m\geq 1$.

\vspace{8pt}
\noindent{\em Proof.}
 The Proposition follows from Chiodo's formula \cite{Chiodo2}
for the Chern characters of $R^*\pi_*\cL$,
\begin{multline} \label{Eq:ChiodosFormula}
{\rm ch}_m(r,k;a_1,\dots,a_n)  =		
		\frac{B_{m+1}(\frac kr)}{(m+1)!} \kappa_m
		- \sum_{i=1}^n 
		\frac{B_{m+1}(\frac{a_i}r)}{(m+1)!} \psi_i^m
		\\  
		+ \frac{r}2 \sum_{w=0}^{r-1} 
		\frac{B_{m+1}(\frac wr)}{(m+1)!} {\xi_{w*}} 
		\left[\frac{(\psi')^m - (-\psi'')^m}{\psi'+\psi''}\right]\, ,
\end{multline}
and from the universal relation 
$$
c(-E^\bullet) = \exp \left(\sum_{m \geq 1} (-1)^m (m-1)! {\rm ch}_m(E^\bullet)\right)\, , 
\ \ \ E^\bullet\in {\mathsf{D}}^b\, .
$$

The factor $r^{|E(\Gamma)|}$ in the Proposition
is due to the factor of $r$ in the last term of Chiodo's formula.{\footnote{The
factor in the last term is $\frac{r}{2}$, but the $\frac{1}{2}$ accounts
for the degree 2 of the boundary map $\xi$ obtained from the ordering of the nodes.}} 
When two such terms corresponding to the same edge are multiplied,
 the factor $r^2$ partly cancels with the factor $\frac{1}{r}$ in the first Chern class of the normal bundle $-\frac{1}{r}(\psi'+ \psi'')$. Thus only a single
 factor of $r$ per edge remains. \qed

\begin{corollary} \label{Cor:ChiodoExp}
The push-forward $\epsilon_* c(-R^*\pi_*L)\in R^*(\oM_{g,n})$ is equal to 
\begin{multline*}
\hspace{-10pt}\sum_{\Gamma\in \mathsf{G}_{g,n}} 
\sum_{w\in \mathsf{W}_{\Gamma,r,k}}
\frac{r^{2g-1-h^1(\Gamma)}}{|\Aut(\Gamma)| }
\;
\xi_{\Gamma*}\Bigg[ \prod_{v \in \V(\Gamma)} e^{-\sum\limits_{m\geq 1} (-1)^{m-1}\frac{B_{m+1}(k/r)}{m(m+1)}\kappa_m(v)} \; \cdot 
\\
\prod_{i=1}^n e^{\sum\limits_{m\geq 1}(-1)^{m-1} \frac{B_{m+1}(a_i/r)}{m(m+1)} \psi^m_{h_i}} \cdot 
\prod_{\substack{e\in \E(\Gamma) \\ e = (h,h')}}
\frac{1-e^{\sum\limits_{m \geq 1} (-1)^{m-1} \frac{B_{m+1}(w(h)/r)}{m(m+1)} [(\psi_h)^m-(-\psi_{h'})^m]}}{\psi_h + \psi_{h'}} \Bigg]\, .
\end{multline*} 
\end{corollary}

\vspace{8pt}
\noindent{\it Proof.} Because the maps $\phi$ and $\varphi$ of Eq.~\eqref{Eq:phi} and~\eqref{Eq:varphi} have degree~1, the degree of the map
\begin{equation} \label{Eq:strataprojection}
\epsilon: \oM_{g;a_1, \dots, a_n}^{r,k} \to \oM_{g,n}
\end{equation}
restricted to the stratum associated to $(\Gamma,w)$ is equal to the product of degrees over vertices, that is,  
$$r^{\sum_{v\in \mathsf{V}(\Gamma)} (2g_v-1)}\, .$$  
The formula of the Corollary is obtained by push-forward from the formula 
of Proposition \ref{Prop:ChiodoExp}. The power of $r$ is given by
\begin{eqnarray*}
|E| + \sum_{v\in \mathsf{V}(\Gamma)} (2g_v-1)& = & |E| - |V| + 
2 \sum_{v\in \mathsf{V}(\Gamma)} g_v\\
& = & (h^1(\Gamma) -1) + 2(g- h^1(\Gamma))\\
& = &2g-1 - h^1(\Gamma)\, .
\end{eqnarray*}

\subsection{Reduction modulo $r$}
\label{Ssec:modr}
Let $g$ be the genus, and let $k\in \Z$ and 
$a_1, \dots, a_n\in \Z$ satisfy
 $$k(2g-2+n)-\sum a_i = 0  \mod r$$ 
for an integer $r>0$. Let 
$$A=(a_1,\ldots, a_n)$$
be the associated vector of $k$-twisted double ramification data.

\vspace{8pt}
\noindent {\bf{Proposition $\bf{3''}$} (Pixton \cite{PixDR2})} {\em 
Let $\Gamma\in\mathsf{G}_{g,n}$ be a fixed stable graph.
Let $\ff$ be a polynomial with $\mathbb{Q}$-coefficients
in variables corresponding to the 
set of half-edges $\mathsf{H}(\Gamma)$. Then the sum 
$$
\mathsf{F}(r) = \sum_{w\in \mathsf{W}_{\Gamma,r,k}} \ff\big(w(h_1), \dots, w(h_{|H(\Gamma)|})\big)
$$
over all possible $k$-weightings mod r 
is (for all sufficiently large $r$)  a {\em polynomial} in $r$ divisible by $r^{h^1(\Gamma)}$.}

\vspace{8pt}

Proposition $3''$ implies Proposition $3'$ (which implies Proposition 3
under the specialization $k=0$). Proposition $3''$ is the basic statement.
See the Appendix for the proof.

%

We will require Proposition $3''$ to prove the following result.
Let $$\epsilon : \oM_{g;a_1, \dots, a_n}^{r,k} \to \oM_{g,n}$$
be the map to the moduli of curves.
 Let $d \geq 0$ be a degree and let $c_d$ denote the $d$-th Chern class.

\begin{prop} \label{Prop:samefreeterm}
The two cycle classes 
$$r^{2d-2g+1} \epsilon_* c_d(-R^*\pi_*\cL)\in R^d(\oM_{g,n})
\ \ \ \text{and} \ \ \ 2^{-d} \, \mathsf{P}_g^{d,r,k}(A)\in R^d(\oM_{g,n})$$
 are polynomials in $r$ (for $r$ sufficiently large).
Moreover,
the two polynomials have the {\em same constant term}.
\end{prop}

\noindent{\em Proof.}
Choose a stable graph $\Gamma$ decorated with  
$\kappa$ classes on the vertices and
$\psi$ classes on half-edges, 
$$[\Gamma,\gamma]\, ,\ \  \ \ \ \gamma=\prod_{v\in \mathsf{V}(\Gamma)} \mathsf{M}_v(\kappa(v))\cdot \prod_{h\in \mathsf{H}(\Gamma)} \psi_h^{m_h}\ .$$
Here, $\mathsf{M}_v$ is a monomial in the $\kappa$ classes at $v$ and $m_h$ is
a non-negative integer for each $h$. The decorated graph data $[\Gamma,\gamma]$  
represents a tautological class of degree $d$. 
We will study the coefficient of the decorated graph in 
the formula for $\epsilon_* c_d(-R^*\pi_*\cL)$ given by 
 Corollary \ref{Cor:ChiodoExp} for $r$ sufficiently large. 

According to Proposition $3''$, the coefficient of $[\Gamma,\gamma]$ 
in the formula for the push-forward $\epsilon_* c_d(-R^*\pi_*\cL)$
is a Laurent polynomial in~$r$. Indeed, the formula
of Corollary \ref{Cor:ChiodoExp} is a combination of a finite number of sums 
of the type of  
Proposition $3''$ multiplied by positive and negative powers of $r$. 
Each contributing sum is obtained by expanding the exponentials in
the formula of Corollary \ref{Cor:ChiodoExp} and selecting a
degree $d$ monomial.

To begin with, let us determine the terms of lowest power in~$r$ in the Laurent polynomial obtained by summation in the formula
of Corollary \eqref{Cor:ChiodoExp}. 
The $(m+1)$st Bernoulli polynomial $B_{m+1}$ has degree~$m+1$. 
Therefore, in the formula of  Corollary \ref{Cor:ChiodoExp},
 {\em every} class of degree $m$, be it $\kappa_m$, $\psi_{h_i}^m$ or an edge{\footnote{The edge contributes 1 to the degree of the class.}} with a class $\psi^{m-1}$, appears with a power of $r$ equal to $\frac{1}{r^{m+1}}$ or 
larger.{\footnote{The power of $\frac{1}{r^{m+1}}$ is $-m-1$. By
{\em larger} here, we mean larger power (so smaller pole or polynomial terms).}} 
If we take a product of $q$ coefficients of expanded exponentials of 
Bernoulli polynomials,
 we will obtain a power of $r$ equal to $\frac{1}{r^{d+q}}$ or larger. 
Thus the lowest possible power of $r$ is obtained if we take a product of $d$ classes of degree 1 each. Then, the power of $r$ is equal to $\frac{1}{r^{2d}}$
which is the lowest possible value. 

After we perform the summation over all $k$-weightings mod $r$ of the graph $\Gamma$, 
we obtain a polynomial in~$r$ divisible by $r^{h^1(\Gamma)}$ by Proposition $3''$. 
The lowest possible power of $r$ 
after summation becomes $r^{h^1(\Gamma)-2d}$. 
Finally, formula of Corollary \ref{Cor:ChiodoExp} carries 
a global factor of $r^{2g-1-h^1(\Gamma)}$. We conclude  $\epsilon_*c_d(-R^*\pi_*\cL)$ is a Laurent polynomial in~$r$ with lowest power of $r$ equal to $r^{2g-2d-1}$. In other words, the product 
$$
r^{2d-2g+1} \epsilon_*c_d(-R^*\pi_*\cL)
$$
is a polynomial in~$r$.

We now identify more precisely the terms of lowest power in $r$.
 As we have seen above, the lowest possible power of $r$ in 
the formula of Corollary \ref{Cor:ChiodoExp} is 
obtained by choosing 
only the terms that contain a degree 1 class multiplied by the second Bernoulli polynomial,
$$
B_2\left(\frac wr\right) = \frac{w^2}{r^2} -\frac w r + \frac16. 
$$
Moreover, in the Bernoulli polynomial only the first term 
$\frac{w^2}{r^2}$ will contribute to the lowest power of~$r$. 
After dropping all the terms which do not contribute to the lowest power of~$r$, the formula of Corollary \ref{Cor:ChiodoExp} simplifies to 
\begin{multline*}
\hspace{-10pt}\sum_{\Gamma\in \mathsf{G}_{g,n}} 
\sum_{w\in \mathsf{W}_{\Gamma,r,k}}
\frac{r^{2g-1-h^1(\Gamma)}}{|\Aut(\Gamma)| }
\;
\xi_{\Gamma*}\Bigg[ \prod_{v \in \V(\Gamma)} \exp \left[-\frac{k^2}{2 r^2}\kappa_1(v)\right]
\prod_{i=1}^n \exp \left[\frac{a_i^2}{2r^2} \psi_{h_i} \right] \cdot 
\\ \hspace{+10pt}
\prod_{e=(h,h')\in \E(\Gamma)}
\frac{1-\exp \left[\frac{w(h)^2}{2r^2} (\psi_h+\psi_{h'})\right]}{\psi_h + \psi_{h'}} \Bigg]\, .
\end{multline*} 
We transform  $w(h)^2$ in the last factor to $-w(h) w(h')$ for symmetry. 
Indeed, 
$$w(h)^2= w(h)(r-w(h')) = - w(h)w(h') \ \mod r\, , $$
so the lowest degree in $r$ is not affected. After 
factoring out the $2r^2$ in the denominators, we  obtain
\begin{multline*}
\hspace{-10pt} r^{2g-1-2d} \cdot 2^{-d} \cdot
\sum_{\Gamma\in \mathsf{G}_{g,n}} 
\sum_{w\in \mathsf{W}_{\Gamma,r,k}}
\frac{r^{-h^1(\Gamma)}}{|\Aut(\Gamma)| }
\;
\xi_{\Gamma*}\Bigg[ \prod_{v \in \V(\Gamma)} e^{-k^2 \kappa_1(v)}
\prod_{i=1}^n e^{a_i^2 \psi_{h_i}} \cdot 
\\ \hspace{+10pt}
\prod_{e=(h,h')\in \E(\Gamma)}
\frac{1-e^{-w(h)w(h’) (\psi_h+\psi_{h'})}}{\psi_h + \psi_{h'}} \Bigg]\, .
\end{multline*} 
Multiplication by $r^{2d-2g+1}$ then yields
exactly $2^{-d} \, \P_g^{d,k,r}(A)$ as claimed. \qed

\section{Localization analysis}

\subsection{Overview}
Let 
$A=(a_1,\ldots,a_n)$ be 
a vector of double ramification data as defined in Section \ref{zz11},
$$\sum_{i=1}^n a_i=0\, .$$
From the positive and negative parts of $A$, we obtain two 
partitions $\mu$ and $\nu$ of the same size
$|\mu|=|\nu|$.
 The double ramification cycle
$$\DR_{g}(A) \in R^g(\oM_{g,n})$$
is defined via the moduli space of stable maps to rubber
$\oM_{g,I}(\PP^1,\mu,\nu)^\sim$
 where $I$ corresponds to the $0$ parts of $A$.
We prove here the claim of Theorem 1,
 $$\DR_{g}(A) = 2^{-g}\, \mathsf{P}_g^g(A) \in R^g(\oM_{g,n})\, .$$
 

\subsection{Target geometry}
Following the notation of \cite{JPT}, let $\PP^1[r]$  
be the projective line with an orbifold 
point $B\Z_r$ at $0\in \PP^1$. 
Let
$$\oM_{g, I, \mu}(\PP^1[r], \nu)$$
be the moduli space 
of stable maps to the orbifold/relative 
pair $(\PP^1[r],\infty)$. The
moduli space 
parameterizes connected, semistable, twisted curves $C$ of genus $g$
with $I$ markings together with a map 
$$f: C \to P\,$$ 
where $P$ is a  destabilization of 
$\PP^1[r]$ over $\infty\in \PP^1[r]$.

We refer the reader to \cite{AGV, JPT} for
the definitions of the moduli space of stable maps to $(\PP^1[r],\infty)$.
The following conditions are required to hold over $0,\infty\in \PP[r]$:
\begin{enumerate}
\item[$\bullet$] The stack structure of the domain curve $C$ occurs
only  at the nodes over $0\in \PP^1[r]$
  and the markings corresponding to $\mu$ (which must to be mapped to
  $0\in \PP^1[r]$). The monodromies associated to the latter
markings are specified by the parts $\mu_i$ of $\mu$. For each part $\mu_i$,
let 
$$\mu_i= \overline{\mu}_i  \mod r  \ \ \ \ \text{where} \ \ 0 \leq \overline{\mu}_i \leq r-1\, .$$

\item[$\bullet$] The map $f$ is finite over $\infty\in \PP^1[r]$ on the last component of $P$ with
  ramification data given by $\nu$. The $\ell(\nu)$ ramification points
are marked.
The map $f$ satisfies the
  ramification matching condition over the internal nodes of the destabilization
$P$.
\end{enumerate}

The moduli space $\oM_{g, I, \mu}(\PP^1[r], \nu)$
has a perfect obstruction theory and a virtual
class of dimension
\begin{equation}\label{zzzz1}
\dim_{\C}\, [\oM_{g, I, \mu}(\PP^1[r], \nu)]^{\vir}= 2g-2+n+\frac{|\nu|}{r}-\sum_{i=1}^{\ell(\mu)} \frac{\overline{\mu}_i}{r}\, ,
\end{equation}
see \cite[Section 1.1]{JPT}.
Recall $n$ is the length of $A$, 
$$n= \ell(\mu)+ \ell(\nu)+\ell(I)\, .$$
Since $|\mu|=|\nu|$, the virtual dimension is an integer. 

We will be most interested in the  case where $r> |\mu|=|\nu|$. Then, 
$$\mu_i=\overline{\mu}_i$$ for all the parts of $\mu$, and 
formula \eqref{zzzz1} for the virtual dimension simplifies to
$$\dim_{\C}\, [\oM_{g, I, \mu}(\PP^1[r], \nu)]^{\vir}= 2g-2+n \, .$$

\subsection{$\C^*$-fixed loci}
The standard $\C^*$-action on $\PP^1$ defined by
$$\xi\cdot[z_0,z_1] =[z_0,\xi z_1]$$
lifts canonically to $\C^*$-actions on
$\PP^1[r]$ and $\oM_{g,I,\mu}(\PP^1[r],\nu)$.
%
%
We describe here the $\C^*$-fixed loci of 
 $\oM_{g,I,\mu}(\PP^1[r],\nu)$.

The $\C^*$-fixed loci of $\oM_{g,I,\mu}(\PP^1[r],\nu)$
are labeled by decorated graphs $\Phi$. 
The vertices $v\in \mathsf{V}(\Phi)$ are
decorated with a genus $\g(v)$ and the legs are labeled with
markings in $\mu \cup I \cup \nu$. Each vertex $v$ is labeled by 
$$0\in\PP^1[r]\ \ \text{or} \ \ \infty\in \PP^1[r]\, .$$
Each edge $e\in \mathsf{E}(\Phi)$ is decorated with a degree $d_e$ that corresponds to the $d_e$-th power map $$\PP^1[r] \to \PP^1[r]\, .$$
The vertex labeling endows the graph $\Phi$ with a {\em bipartite} structure.

A vertex of $\Phi$ over $0\in \PP^1[r]$ corresponds to a stable map 
contracted to $0$ given by an element of $\oM_{\g(v), I(v),\mu(v)}(B\Z_r)$
with monodromies $\mu(v)$ specified by the corresponding entries of 
$\mu$ {\em and} the negatives of the degrees $d_e$ of the incident edges 
(modulo $r$). 
We will use the notation
$$\oM_v^{r}\, =\, \oM^r_{\g(v);I(v),\mu(v)}\, =\, \oM_{\g(v),I(v),\mu(v)}(B\Z_r)\, .$$
The fundamental class  
 $\left[ \oM_v^{r} \right]$ of $\oM_v^r$ is of dimension
$$\dim_{\C}\, \left[ \oM_v^{r} \right]=
3\g(v)-3+ \ell(I(v))+\ell(\mu(v))\,  .$$
 As  explained in Section \ref{cmc}, 
 $$\epsilon:\oM_v^r\rightarrow \oM_{\g(v),\n(v)}$$
is a finite covering of the moduli space $\oM_{\g(v), \n(v)}$ 
where $$\n(v)=\ell(\mu(v))+\ell(I(v))$$ is 
the total number of markings and edges
incident to $v$. 

The stable map over $\infty\in \PP^1[r]$ can take two different forms. If the target does not expand then the stable map has $\ell(\nu)$ preimages of $\infty \in \PP^1[r]$. Each preimage is described by an unstable vertex of $\Phi$ labeled by~$\infty$. If the target expands then the stable map is a  possibly disconnected rubber map. The ramification data is given by incident edges over 0 of the rubber and by elements of $\nu$ over $\infty$ of the rubber. In this case every vertex $v$ over $\infty\in \PP^1[r]$ describes a connected component of the rubber map. 

Let $\g(\infty)$ be the genus of the possibly disconnected domain of the rubber map. We will denote{\footnote{We omit the data $\g(\infty)$, $I(\infty)$, $\nu$, and $\delta(\infty)$
in the notation.}} the moduli space
of stable maps to rubber by $\oM_\infty^{\sim}$ and the
 virtual fundamental class by $\left[\oM_\infty^{\sim} \right]^{\vir}$,
$$\dim_{\C}\, \left[ \oM_\infty^{\sim} \right]^{\vir}=
2\g(\infty)-3+ \n(\infty)\,  ,$$
where $\n(\infty)$ is the total number of markings and edges
incident to vertices over~$\infty$,
$$\n(\infty)= \ell(I(\infty))+ \ell(\nu)+|\mathsf{E}(\Phi)|\, .$$ 
The image of the virtual fundamental class 
$\left[\oM_\infty^{\sim} \right]^{\vir}$
in the moduli space of (not necessarily connected{\footnote{By \cite[Lemma 2.3]{BSSZ}, a double ramification cycle
 vanishes as soon as there are at least two nontrivial components
of the domain. However, we will {\em not} require the
vanishing result here.
Our analysis will naturally avoid
disconnected domains mapping to the rubber over $\infty\in \PP^1[r]$.}}) stable curves is denoted by $\DR_\infty$.

A vertex $v\in\mathsf{V}(\Phi)$ is {\em unstable} when $2 \g(v) - 2 + \n(v) \leq 0$. There are four types of unstable vertices: 
\begin{enumerate}
\item[(i)] $v\mapsto 0$, $\g(v)=0$, $v$ carries no markings and
 one incident edge,
\item[(ii)] $v\mapsto 0$, $\g(v)=0$,  
$v$ carries no markings and two incident edges,
\item[(iii)] $v\mapsto 0$, $\g(v)=0$, 
 $v$ carries one marking and one incident edge,
\item[(iv)] $v \mapsto \infty$, $\g(v)=0$, 
$v$ carries one marking and one incident edge.
\end{enumerate}
The target of the stable map expands iff there is at least one stable vertex over~$\infty$.

A stable map in the $\C^*$-fixed locus
corresponding to  $\Phi$ is obtained by gluing together maps
associated to the vertices $v\in \mathsf{V}(\Phi)$ with
Galois covers 
associated to the edges. Denote by $V_{\rm st}^0(\Phi)$ the set of stable vertices of $\Phi$ over~0. Then the $\C^*$-fixed locus
corresponding to $\Phi$ is isomorphic to the product 
$$
\oM_\Phi =
 \begin{cases}
      \prod\limits_{v\in \mathsf{V}^0_{\rm st}(\Phi)}\, \oM_v^r\ \times\ \oM_\infty^{\sim}      \, , & \text{if the target expands,} 
\vspace{7pt}\\ 
    \prod\limits_{v\in \mathsf{V}^0_{\rm st}(\Phi)}\, \oM_v^{r}
\, , & \text{if the target does not expand,}
  \end{cases}
$$
quotiented by the automorphism group of $\Phi$ {\em and} the product of 
cyclic groups $\Z_{d_e}$ associated to the Galois covers of the edges.

Thus the natural morphism corresponding to $\Phi$,
$$\iota: \oM_\Phi \rightarrow \oM_{g, I, \mu}(\PP^1[r], \nu)\, ,$$
is of degree 
$$
|\Aut(\Phi)| \prod_{e\in \mathsf{E}(\Phi)}d_e.
$$
onto the image $\iota(\oM_{\Phi})$.


\vspace{5pt}

Recall that $D = |\mu| = |\nu|$ is the degree of the stable map.

\begin{lemma}
For $r > D$, the unstable vertices of type (i) and (ii)
can {\em not} occur.
\end{lemma}

\vspace{6pt}
\noindent{\em Proof.} At each stable vertex $v\in \mathsf{V}(\Phi)$
over $0\in \PP^1[r]$, the condition
\begin{equation*}
  |\mu_v| = 0   \mod {r}
\end{equation*}
has to hold in order for $\oM_{\g_v, I_v, \mu_v}(B\Z_r)$ to be
non-empty. Because of the balancing condition at the nodes, the parts
of $\mu$ at $v$ have to add up to the sum of the degrees of the
incident edges modulo $r$. Since $r > D$, the parts of $\mu$ at $v$
must add up to {\em exactly} the sum of the degrees of the incident edges.
So the incident edges of the stable vertices and
the type (iii) unstable vertices must account for the 
{\em entire} degree $D$. Thus, there is no remaining
edge degree for unstable vertices of type (i) and (ii). \qed

\subsection{Localization formula}
\label{locformula}

We write the $\C^*$-equivariant Chow ring of a point
as $$A^*_{\C^*}(\bullet) = \Q[t]\, ,$$
where $t$ is the first Chern class of the standard representation.

For the localization formula, we will
require the inverse of the $\C^*$-equivariant Euler
class of the virtual normal bundle in $\oM_{g,I,\mu}(\PP^1[r],\nu)$
to the $\C^*$-fixed
locus corresponding  $\Phi$. Let
$$f: C \to \PP[r]\, , \ \ \ \ \ [f]\in \oM_\Phi\, .$$
The inverse Euler class is represented by 
\begin{equation}
  \label{eq:contribs}
\frac{1}{e(\text{Norm}^{\vir})}=  \frac{e(H^1(C, f^*T_{\PP[r]} (-\infty)))}{e(H^0(C, f^*T_{\PP[r]} (-\infty)))} \frac 1{\prod_i e(N_i)} \frac 1{e(N_\infty)}\, .
\end{equation}

Formula \eqref{eq:contribs} has several terms which require
explanation. We assume $$r>|\mu|=|\nu|\, ,$$ so unstable vertices over
$0\in \PP^1[r]$ of type (i) and (ii) do not occur in $\Phi$. 
First consider  the leading factor of~\eqref{eq:contribs},
\begin{equation}\label{xhh}
\frac{e(H^1(C, f^*T_{\PP[r]} (-\infty)))}{e(H^0(C, f^*T_{\PP[r]} (-\infty)))}\, .
\end{equation}
To compute this factor one uses the normalization exact sequence for the domain $C$ tensored with the line bundle $f^*T_{\PP[r]}(-\infty)$. The associated long exact sequence in cohomology decomposes the leading factor into a product of vertex,
edge, and node contributions.

\begin{enumerate}

\item[$\bullet$]
Let $v\in \mathsf{V}(\Phi)$ be a stable vertex over $0\in\PP^1[r]$ 
corresponding to 
a contracted component $C_v$.
From the resulting stable map to the
orbifold point 
$$C_v \rightarrow B\Z_r \subset \PP[r]\, , $$
we obtain an orbifold line bundle
$$L=f^*T_{\PP^1[r]}|_0$$ 
on $C_v$ which is an $r$th root of $\cO_{C_v}$. 
Therefore,  the
contribution 
\begin{equation*}
  \frac{e(H^1(C_v, f^*T_{\PP[r]}|_0))}{e(H^0(C_v, f^*T_{\PP[r]}|_0))}
\end{equation*}
yields the class
\begin{equation*}  
c_{\mathrm{rk}}\left((-R^*\pi_*\cL) \otimes \cO^{(1/r)}\right) \
\in A^*(\oM^r_{\g(v),I(v),\mu(v)}) \otimes \Q\left[t,\frac 1t\right]\, ,
\end{equation*}
where $\cL\rightarrow \oM^r_{\g(v),I(v),\mu(v)}$ 
is the universal $r$th root, $\cO^{(1/r)}$ is a trivial
line bundle with a $\C^*$-action of weight $\frac 1r$, and
$$\mathrm{rk} = g_v-1 + |\mathsf{E}(v)|$$
is the virtual rank of $-R^*\pi_*\cL$.

Unstable vertices of type (iii) over $0\in \PP^1[r]$ and
the vertices over $\infty$ contribute factors of 1.

\item[$\bullet$] The edge contribution is trivial since the degree $\frac {d_i}r$ of
$f^*T_{\PP[r]} (-\infty)$ is less than 1, see \cite[Section 2.2]{JPT}.

\item[$\bullet$] The contribution of a node $N$ over $0\in \PP^1[r]$ is trivial. The
space of sections
 $H^0(N,f^*T_{\PP[r]} (-\infty))$ vanishes  because $N$ must be
stacky, and $H^1(N,f^*T_{\PP[r]} (-\infty))$ is trivial for
dimension reasons.

Nodes over $\infty\in \PP^1[r]$ contribute 1.
\end{enumerate}

\noindent Consider next the last two factors of \eqref{eq:contribs},
$$\frac 1{\prod_i e(N_i)} \frac 1{e(N_\infty)}\, .$$

\begin{enumerate}
\item[$\bullet$]
The  product ${\prod_i e(N_i)^{-1}}$
is over all nodes over $0\in \PP^1[r]$ of the domain $C$ 
which are {\em forced} by the graph $\Phi$.
These are nodes where the edges of $\Phi$ are attached to the vertices.
If $N$ is a node over $0$ forced by the edge $e\in \mathsf{E}(\Phi)$ 
and the associated vertex~$v$ is stable, then 
\begin{equation} \label{Eq:eN}
 e(N) = \frac t{r \, d_e} - \frac{\psi_e}r.
\end{equation}
This expression corresponds
 to the smoothing of the node $N$ of the domain curve: $e(N)$ is 
the first Chern classes of the normal bundle of 
the divisors of nodal domain curves.
The first Chern classes of the cotangent lines to the branches at the node 
are divided by $r$ because of the orbifold twist, see Section
\ref{Ssec:koverr}.

In the case of an unstable vertex of type (iii), the associated
edge does {\em not}
produce a node of the domain. The type (iii) edge incidences do
{\em not} appear in ${\prod_i e(N_i)^{-1}}$.

\item[$\bullet$]
$N_\infty$ corresponds to the target degeneration over $\infty \in \PP^1[r]$. The factor $e(N_\infty)$
is 1 if the {\em target} $(\PP^1[r], \infty)$ does not degenerate and 
$$e(N_\infty)= -\frac{t + \psi_\infty}{{\prod_{e\in\mathsf{E}(\Phi)}d_e}}$$
if the target does degenerate \cite{GrV}. Here, $\psi_\infty$ is the first Chern class of the cotangent line bundle to the bubble in the target curve at the attachment point to $\infty$ of $\PP[r]$. 
\end{enumerate}

Summing up our analysis, we write out the outcome of the virtual localization formula~\cite{GrP} in $A^*(\oM_{g, I, \mu}(\PP^1[r], \nu))\otimes  \Q[t,\frac{1}{t}]$. We have
\begin{equation}\label{llff}
[\oM_{g, I, \mu}(\PP^1[r], \nu)]^{\vir} =
\sum_{\Phi} \frac{1}{|\Aut(\Phi)|} \, \frac{1}
{\prod_{e\in\mathsf{E}(\Phi)}d_e}\cdot 
\iota_*\left(\frac{[\oM_\Phi]^{\vir}}
{e(\text{Norm}^{\vir})}\right)\,,
\end{equation}
where $\frac{[\oM_\Phi]^{\vir}}
{e(\text{Norm}^{\vir})}$ is the product of the following factors:
\begin{itemize}
\item 
$ \displaystyle
\prod_{e \in \mathsf{E}(v)} \frac {r}{\frac{t}{d_e}-\psi_e}\, \cdot\, 
\sum_{d \geq 0} c_d(-R^*\pi_*\cL) \left(\frac{t}r\right)^{\g(v)-1+|\mathsf{E}(v)|-d} \;
$
\begin{minipage}{10em}
for each stable vertex $v\in \mathsf{V}(\Phi)$ over 0,
\end{minipage}
\item 
$\displaystyle
   -\frac {{\prod_{e\in\mathsf{E}(\Phi)}d_e}}{t+ \psi_\infty} 
$
\hspace{16em} 
\begin{minipage}{10em}
if the target degenerates.
\end{minipage}

\end{itemize}
\begin{remark} \label{Rem:2conventionsB}
As we mentioned in Remark~\ref{Rem:2conventionsA}, there are two conventions for the orbifold structure on the space of $r$th roots. In the paper,
we follow the convention which assigns to each node of the curve the stabilizer $\Z/r\Z$. It is also possible to assign to a node of type $(a',a'')$ the stabilizer $\Z/q\Z$, where $q$ is the order of $a'$ and $a''$ in $\Z/r \Z$. 
We list here the required modifications of our formulas 
when the latter convention is followed. 

Denote $p = \gcd(r,a') = \gcd(r,a'')$. We then have $pq=r$. In Section~\ref{Sssec:normalbundle}, the first Chern class of the normal line bundle to a boundary divisor in the space of $r$th roots becomes $-(\psi'+\psi'')/q$ rather than $-(\psi'+\psi'')/r$. The degrees of the maps $\phi$ and $\varphi$ in Equations~\eqref{Eq:phi} and \eqref{Eq:varphi} is $p$ instead of~1. The factor $\frac{r}2$ in the last term of Chiodo's formula~\eqref{Eq:ChiodosFormula} becomes $\frac{q}2$. Because of the above, the factor $r^{|E(\Gamma)|}$ in Proposition~\ref{Prop:ChiodoExp} is transformed into $q^{|E(\Gamma)|}$. The statement of Corollary~\ref{Cor:ChiodoExp} is unchanged, but in the proof the factor $r^{|E(\Gamma)|}$ is now obtained as a product of $q^{|E(\Gamma)|}$ from Proposition~\ref{Prop:ChiodoExp} and $p^{|E(\Gamma)|}$ from the degree of $\phi$ and $\varphi$. 

In the localization formula, the Euler classes of the normal bundles at the nodes over~0 become $\frac{t}{q d_e} - \frac{\psi_e}q$ instead of $\frac{t}{r d_e} - \frac{\psi_e}r$. On the other hand, at each node there are $p$ ways to reconstruct an $r$th root bundle from its restrictions to the two branches. Thus the contribution of each node $\frac{p}{e(N)}$ still has an $r$ in the numerator.
\end{remark}

\subsection{Extracting the double ramification cycle}

\subsubsection{Three operations}
We will now perform three operations on the localization formula \eqref{llff}
for the virtual class $[\oM_{g, I, \mu}(\PP^1[r], \nu)]^{\vir}$:

\begin{enumerate}
\item[(i)] the $\C^*$-equivariant push-forward via
\begin{equation}\label{xx88}
\epsilon: \oM_{g,I,\mu}(\PP^1[r],\nu) \rightarrow \oM_{g,n}
\end{equation}
to the moduli space $\oM_{g,n}$ with trivial $\C^*$-action, 
\item[(ii)] extraction of the coefficient of $t^{-1}$ after push-forward by $\epsilon_*$, 
\item[(iii)] extraction of the coefficient of $r^0$.
\end{enumerate}

After push-forward by $\epsilon_*$, the coefficient of $t^{-1}$
is equal to 0 because 
$$\epsilon_*[\oM_{g, I, \mu}(\PP^1[r], \nu)]^{\vir} \in A^*(\oM_{g,n})\otimes
\Q[t]\, .$$
Using the result of Section \ref{Ssec:modr}, all terms of
the $t^{-1}$ coefficient will be seen to be 
polynomials in~$r$,
so operation~(iii) will be well-defined. 
After operations~(i-iii), only two nonzero terms will remain.
The cancellation of the two remaining terms will prove Theorem \ref{FFFF}.

To perform (i-iii),
we multiply the $\epsilon$-push-forward of the
localization formula \eqref{llff} by $t$
 and extract the coefficient of $t^0 r^0$.
To simplify the computations, we introduce the new variable 
$$s=tr\, .$$ 
Then, instead of extracting the coefficient of $t^0r^0$,
we extract the coefficient of $s^0 r^0$. 

\subsubsection{Push-forward to $\oM_{g,n}$} \label{pppfff}

For each vertex $v\in \mathsf{V}(\Phi)$, we have
$$\oM_v^r =\oM^r_{\g(v);I(v),\mu(v)}\, .$$
Following  Section \ref{cmc}, we denote
the map to moduli by
\begin{equation} \label{xx99}
\epsilon: \oM^r_{\g(v);I(v),\mu(v)} \rightarrow \oM_{\g(v),\n(v)}\, .\end{equation}
The use of the notation $\epsilon$ in \eqref{xx88} and \eqref{xx99}
is compatible. Denote by 
\begin{equation}\label{ffrrtt}
\widehat{c}_d \, = \, r^{2d-2\g(v)+1} \epsilon_* c_d(-R^*\pi_*\cL)
\, \in\, R^d(\oM_{\g(v),\n(v)})\, .
\end{equation}
By Proposition \ref{Prop:samefreeterm}, 
$\widehat{c}_d$  is a polynomial in $r$
for $r$ sufficiently large.

From \eqref{llff} and the contribution calculus for $\Phi$ 
presented in Section \ref{pppfff}, we have a complete formula
for the $\C^*$-equivariant push-forward of
$t$ times the virtual class:
\begin{multline}\label{llfff}
\epsilon_*\Big(t[\oM_{g, I, \mu}(\PP^1[r], \nu)]^{\vir}\Big) = \\
\frac{s}{r}\cdot\sum_{\Phi} \frac{1}{|\Aut(\Phi)|} \, \frac{1}
{\prod_{e\in\mathsf{E}(\Phi)}d_e}\cdot 
\epsilon_*\iota_*\left(\frac{[\oM_\Phi]^{\vir}}
{e(\text{Norm}^{\vir})}\right)\, ,
\end{multline}
where $\epsilon_*\iota_*\frac{[\oM_\Phi]^{\vir}}
{e(\text{Norm}^{\vir})}$ is the product of the following factors:
\begin{itemize}
\item a factor
\begin{equation*}
  \frac{r}{s}\cdot \prod_{e \in \mathsf{E}(v)} \frac {d_e}{1-\frac{rd_e}{s}\psi_e}\, \cdot\,
  \sum_{d \geq 0} \widehat{c}_d \, s^{\g(v)-d}\ \ \in R^d(\oM_{\g(v),\n(v)})\otimes \Q\left[s,\frac 1s\right]
\end{equation*}
for each stable vertex $v\in \mathsf{V}(\Phi)$ over 0,
\item a factor 
\begin{equation*}
   -\frac{r}{s}\cdot \frac {{\prod_{e\in\mathsf{E}(\Phi)}d_e}       }{1+\frac{r}{s} \psi_\infty}\cdot \mathsf{DR}_\infty 
\end{equation*}
if the target degenerates.
\end{itemize}


\subsubsection{Extracting coefficients}

\paragraph{Extracting the coefficient of $r^0$.} 
By Proposition \ref{Prop:samefreeterm}, 
the classes $\widehat{c}_d$ are polynomial in $r$ for $r$ sufficiently large.
We have an $r$ in the denominator in the prefactor
on the right side of \eqref{llfff} which
 comes from the multiplication by $t$ on the left side.
However, in all other factors, we only have positive powers of $r$,
 with at least one $r$ per stable vertex of the graph over~0 and one more~$r$ if the target degenerates. The {\em only} graphs $\Phi$ which
 contribute to the coefficient of $r^0$ are those with 
{\em exactly one $r$ in the numerator}. There are only two graphs which have exactly one $r$ factor in the numerator:
\begin{enumerate}
\item[$\bullet$]
the graph $\Phi'$ with a stable vertex of full genus~$g$ over~0 and 
$\ell(\nu)$ type~(iv) unstable vertices over~$\infty$, 
\item[$\bullet$]the graph $\Phi''$ with a stable vertex of full genus~$g$ over~$\infty$ and $\ell(\mu)$ type~(iii) unstable vertices over 0. 
\end{enumerate}

No terms involving $\psi$ classes contribute
to the $r^0$ coefficient of either $\Phi'$ or $\Phi''$
since every $\psi$ class in the 
localization formula comes with an extra factor of $r$.
We can now write the $r^0$ coefficient of the right side of \eqref{llfff}:
\begin{equation}\label{llffff}
\text{Coeff}_{r^0}\left[\epsilon_*\Big(t[\oM_{g, I, \mu}(\PP^1[r], \nu)]^{\vir}\Big)\right] \ =\ 
\text{Coeff}_{r^0}\left[\sum_{d \geq 0} \widehat{c}_d \, s^{g-d}\right] \ -\ 
   \mathsf{DR}_g(A)\, 
 \end{equation}
in $R^d(\oM_{g,n})\otimes \Q[s,\frac{1}{s}]$.

\paragraph{Extracting the coefficient of $s^0$.} 
The remaining powers of $s$ in \eqref{llffff} appear only
 with the classes $\widehat{c}_d$ in contribution of the graph $\Phi'$.
In order to obtain $s^0$, we take to take $d=g$,
\begin{equation}\label{llfffff}
\text{Coeff}_{s^0r^0}\left[\epsilon_*\Big(t[\oM_{g, I, \mu}(\PP^1[r], \nu)]^{\vir}\Big)\right] \ =\ 
\text{Coeff}_{r^0}\left[\widehat{c}_g\right]  \ -\ 
   \mathsf{DR}_g(A)\, 
 \end{equation}
in $R^d(\oM_{g,n})$.

\subsubsection{Final relation}
Since $\text{Coeff}_{s^0r^0}\left[\epsilon_*\Big(t[\oM_{g, I, \mu}(\PP^1[r], \nu)]^{\vir}\Big)\right]$ vanishes, we can rewrite \eqref{llfffff} as
\begin{equation}\label{jjqq}
   \mathsf{DR}_g(A) = \, \text{Coeff}_{r^0}\left[r^{2g-2g+1}
\epsilon_* c_g(-R^*\pi_*\cL)\right]
\, \in\, R^g(\oM_{g,n})\, ,
 \end{equation}
using definition \eqref{ffrrtt} of $\widehat{c}_d$. 
By Proposition \ref{Prop:samefreeterm} applied to the right side of relation 
\eqref{jjqq}, 
$$\mathsf{DR}_g(A) =  2^{-g}\, \mathsf{P}_g^g(A)
\, \in\, R^g(\oM_{g,n})\, . $$
The proof of Theorem \ref{FFFF} is complete. \qed

\section{A few applications}

\subsection{A new expression for $\lambda_g$}

We discuss here a few more examples of the formula for 
Corollary~\ref{xzz}: a formula for the class $\lambda_g$ supported on the boundary divisor with a nonseparating node. We give the formulas for the class $\lambda_g$ in $\oM_g$ up to genus~4. 

The answer is expressed as a sum of labeled stable graphs with coefficients where each stable graph has at least one nontrivial cycle. The genus of each vertex is written inside the vertex. The powers of $\psi$-classes are written on the corresponding half-edges (zero powers are omitted). The expressions do not involve $\kappa$-classes. 

Each labeled graph $\Gamma$ describes a moduli space $\oM_\Gamma$ (a product of moduli spaces associated with its vertices), a tautological class $\alpha \in R^*(\oM_\Gamma)$ and a natural map $\pi: \oM_\Gamma \to 
\oM_g$. Our convention is that $\Gamma$ then represents the cycle class $\pi_* \alpha$. For instance, if the graph carries no $\psi$-classes, the class $\alpha$ equals~1 and the map $\pi$ is of degree $|\Aut (\Gamma)|$. The cycle class represented by $\Gamma$ is then $|\Aut (\Gamma)|$ times the class of the image of~$\pi$.

\tikz{\coordinate (A) at (0,0); \coordinate (B) at (1,0); \coordinate (C) at (0.6,0.5);}

\tikzset{baseline=0, label distance=-3mm}
\def\NC{\draw (0,0.25) circle(0.25);}
\def\NL{\draw plot [smooth,tension=1.5] coordinates {(0,0) (-0.2,0.5) (-0.5,0.2) (0,0)};}
\def\NR{\draw plot [smooth,tension=1.5] coordinates {(0,0) (0.2,0.5) (0.5,0.2) (0,0)};}
\def\NN{\NL\NR}
\def\NNN{\NN \begin{scope}[rotate=180] \NR \end{scope}}
\def\NNNN{\NN \begin{scope}[rotate=180] \NN \end{scope}}
\def\NRS{\begin{scope}[shift={(B)}] \NR \end{scope}}
\def\NRD{\begin{scope}[rotate around={-90:(B)}] \NRS \end{scope}}
\def\DE{\draw plot [smooth,tension=1] coordinates {(0,0) (0.5,0.15) (1,0)}; \draw plot [smooth,tension=1.5] coordinates {(0,0) (0.5,-0.15) (1,0)};}
\def\DES{\begin{scope}[shift={(B)}] \DE \end{scope}}
\def\TE{\DE \draw (A) -- (B);}
\def\QE{\DE \draw plot [smooth,tension=1] coordinates {(A) (0.5,0.05) (B)}; \draw plot [smooth,tension=1.5] coordinates {(A) (0.5,-0.05) (B)};}
\def\T{\draw (0.2,0) -- (C) -- (B) -- (0.2,0);}
\def\TT{\draw (C) -- (B) -- (0.2,0); \draw plot [smooth,tension=1] coordinates {(0.2,0) (0.3,0.3) (C)}; \draw plot [smooth,tension=1] coordinates {(0.2,0) (0.5,0.2) (C)};}
\newcommand{\nn}[3]{\draw (#1)++(#2:3mm) node[fill=white,fill opacity=.85,inner sep=0mm,text=black,text opacity=1] {$\substack{\psi^#3}$};}
\renewcommand{\gg}[2]{\fill (#2) circle(1.3mm) node {\color{white}$\substack #1$};}

\paragraph{Genus~1.}
\begin{equation*}
 \lambda_1 = \frac 1{24} \tikz{\NC \gg{0}{A}}.
\end{equation*}

\paragraph{Genus~2.}
\begin{equation*}
\lambda_2 = 
  \frac 1{240} \tikz{\NC \nn{A}{130}{{}} \gg{1}{A}}
  + \frac 1{1152} \tikz{\NN \gg{0}{A}}.
\end{equation*}

\paragraph{Genus~3.}
\begin{align*}
\lambda_3 &= 
  \frac 1{2016} \tikz{\NC \nn{A}{130}{2} \gg{2}{A}}
  + \frac 1{2016} \tikz{\NC \nn{A}{130}{{}} \nn{A}{50}{{}} \gg{3}{A}}
  - \frac 1{672} \tikz{\DE \nn{A}{30}{{}} \gg{1}{A} \gg{1}{B}}
  + \frac 1{5760} \tikz{\NN \nn{A}{160}{{}} \gg{1}{A}} \\
  &
  - \frac{13}{30240} \tikz{\TE \gg{0}{A} \gg{1}{B}}
  - \frac 1{5760} \tikz{\NL \DE \gg{0}{A} \gg{1}{B}}
  + \frac 1{82944} \tikz{\NNN \gg{0}{A}}.
\end{align*}

\paragraph{Genus~4.}
\begin{align*}
\lambda_4 &= 
  \frac 1{11520} \tikz{\NC \nn{A}{130}{3} \gg{3}{A}}
  + \frac 1{3840} \tikz{\NC \nn{A}{130}{2} \nn{A}{50}{{}} \gg{3}{A}}
  - \frac 1{2880} \tikz{\DE \nn{A}{30}{2} \gg{1}{A} \gg{2}{B}}
  - \frac 1{3840} \tikz{\DE \nn{A}{30}{{}} \nn{A}{-30}{{}} \gg{1}{A} \gg{2}{B}}
  - \frac 1{1440} \tikz{\DE \nn{A}{30}{{}} \nn{B}{150}{{}} \gg{1}{A} \gg{2}{B}} \\
  &- \frac 1{1920} \tikz{\DE \nn{A}{30}{{}} \nn{B}{-150}{{}} \gg{1}{A} \gg{2}{B}}
  - \frac 1{2880} \tikz{\DE \nn{B}{150}{2} \gg{1}{A} \gg{2}{B}}
  - \frac 1{3840} \tikz{\DE \nn{B}{150}{{}} \nn{B}{-150}{{}} \gg{1}{A} \gg{2}{B}}
  + \frac 1{48384} \tikz{\NN \nn{A}{160}{2} \gg{2}{A}}
  + \frac 1{48384} \tikz{\NN \nn{A}{160}{{}} \nn{A}{110}{{}} \gg{2}{A}} \\
  &+ \frac 1{115200} \tikz{\NN \nn{A}{160}{{}} \nn{A}{20}{{}} \gg{2}{A}}
  + \frac 1{960} \tikz{\T \nn{B}{180}{{}} \gg{1}{B} \gg{1}{0.2,0} \gg{1}{C}}
  - \frac{23}{100800} \tikz{\TE \nn{A}{-30}{{}} \gg{2}{A} \gg{0}{B}}
  - \frac 1{57600} \tikz{\DE \NRS \nn{B}{20}{{}} \gg{2}{A} \gg{0}{B}} \\
&  - \frac 1{16128} \tikz{\DE \NRS \nn{B}{-150}{{}} \gg{2}{A} \gg{0}{B}}
  - \frac 1{16128} \tikz{\DE \NRS \nn{A}{-30}{{}} \gg{2}{A} \gg{0}{B}}
  - \frac 1{57600} \tikz{\DE \NRS \nn{B}{20}{{}} \gg{1}{A} \gg{1}{B}}
  - \frac 1{16128} \tikz{\DE \NRS \nn{B}{-150}{{}} \gg{1}{A} \gg{1}{B}} \\
  & - \frac 1{16128} \tikz{\DE \NRS \nn{A}{-30}{{}} \gg{1}{A} \gg{1}{B}}
  - \frac{23}{100800} \tikz{\TE \nn{A}{-30}{{}} \gg{1}{A} \gg{1}{B}}
  + \frac{23}{100800} \tikz{\TT \gg{2}{B} \gg{0}{0.2,0} \gg{0}{C}}
  + \frac{23}{50400} \tikz{\TT \gg{1}{B} \gg{1}{0.2,0} \gg{0}{C}}
  + \frac 1{16128} \tikz{\T \NRS \gg{0}{B} \gg{1}{0.2,0} \gg{1}{C}} \\
  &
  + \frac 1{115200} \tikz{\DE \DES \gg{1}{A} \gg{0}{B} \gg{1}{2,0}}
  + \frac 1{276480} \tikz{\NNN \nn{A}{20}{{}} \gg{1}{A}}
  - \frac{13}{725760} \tikz{\NL \TE \gg{1}{A} \gg{0}{B}}
  - \frac 1{138240} \tikz{\NL \NRS \DE \gg{1}{A} \gg{0}{B}} \\
  &
  - \frac{43}{1612800} \tikz{\QE \gg{1}{A} \gg{0}{B}}
  - \frac{13}{725760} \tikz{\NRS \TE \gg{1}{A} \gg{0}{B}}
  - \frac 1{276480} \tikz{\NRS \NRD \DE \gg{1}{A} \gg{0}{B}}
  - \frac 1{7962624} \tikz{\NNNN \gg{0}{A}}
\end{align*}

All these expressions are obtained by substituting $a_1 = \dots = a_n = 0$ into the formula for the double ramification cycle.

\subsection{Hodge integrals}

Corollary~\ref{xzz} may be applied to any Hodge integral which contains
the top Chern class of the Hodge bundle. As an example, we consider the
the following Hodge integral (first calculated in \cite{FP-H}) related to constant map contribution in
the Gromov-Witten of Calabi-Yau threefolds.

\begin{prop}
For $g\geq 1$, we have
$$
\int\limits_{\oM_{g+1}} \lambda_{g+1} \lambda_g \lambda_{g-1} = 
-\frac12 \frac1{(2g)!} \frac{B_{2g}}{2g} \frac{B_{2g+2}}{2g+2}\, .
$$
\end{prop}

\vspace{8pt}
\noindent{\em Proof.}
We compute the integral by replacing $\lambda_{g+1}$ with the expression
obtained from Corollary \ref{xzz}.

In Pixton's formula for $\DR_{g+1}(\emptyset)=(-1)^{g+1}\lambda_{g+1}$,
we only have to consider graphs with no separating edges. 
On the other hand, $\lambda_g \lambda_{g-1}$ vanishes on strata
 as soon as there are at least {\em two} nonseparating edges. 
Thus we are left with just one graph: the graph with a vertex of genus $g$ and  one loop. 

The edge term in Pixton's formula for 
$(-1)^{g+1}\DR_{g+1}(\emptyset)=\lambda_{g+1}$ 
for the unique remaining 
graph is 
$$
(-1)^{g+1}\frac{1-e^{\frac{a^2}2(\psi'+\psi'')}}{\psi'+\psi''}
= 
(-1)^g\sum_{k \geq 0} \frac{a^{2k+2}}{2^{k+1}\; (k+1)!}
(\psi'+ \psi'')^k.
$$
 The contribution of the graph is therefore the $r$-free term of
$$
 \frac{(-1)^g}{2r} \sum_{a=0}^{r-1}  \frac{a^{2g+2}}{2^{g+1}\; (g+1)!}\; \cdot \;  \int\limits_{\oM_{g,2}} (\psi_1+ \psi_2)^g \lambda_g \lambda_{g-1}\, .
$$
The factor of $\frac{1}{2}$ comes from the automorphism group of the graph, and
 the factor of $\frac{1}{r}$ comes from
the first Betti number of graph. 
Since
$$
\sum_{a=0}^{r-1} a^{2g+2} = B_{2g+2} r + O(r^2),  
$$
we see that the $r$-free term equals
$$
(-1)^g \frac{B_{2g+2}}{2^{g+2}\, (g+1)!} \int\limits_{\oM_{g,2}} (\psi_1+ \psi_2)^g \lambda_g \lambda_{g-1}.
$$
Using Lemma~\ref{ggww} below we can rewrite this expression as
$$
(-1)^g \frac{B_{2g+2}}{2^{g+2}\, (g+1)!} \cdot
(-1)^{g+1} \frac{B_{2g}}{2g} \frac1{(2g-1)!!}
=
-\frac12 \frac1{(2g)!} \frac{B_{2g}}{2g} \frac{B_{2g+2}}{2g+2}\, ,
$$
which completes the derivation. \qed


\begin{lemma} \label{ggww}
For $g\geq 1 $, we have
$$
\int\limits_{\oM_{g,2}} (\psi_1 + \psi_2)^g \lambda_g \lambda_{g-1} = 
(-1)^{g+1} \frac{B_{2g}}{2g} \frac1{(2g-1)!!}\, .
$$
\end{lemma}

\vspace{8pt}
\noindent{\em Proof.}
According to Faber's socle formula \cite{Fab,GetP}, we have
\begin{equation} \label{Eq:socle}
\int\limits_{\oM_{g,2}} \psi_1^p \psi_2^q \lambda_g \lambda_{g-1} 
= \frac{(-1)^{g+1} \, B_{2g}}{2^{2g}\, g \, (2p-1)!! (2q-1)!!}
\end{equation}
whenever $p+q = g$. The formula holds even if $p$ or $q$ is equal to~0.{\footnote{The socle formula is true if at most one power of a $\psi$ class vanishes 
while the powers of all other $\psi$ classes are positive.}} 
After expanding $(\psi_1+\psi_2)^g$ and applying \eqref{Eq:socle},
 we obtain 
$$
\int\limits_{\oM_{g,2}} (\psi_1 + \psi_2)^g \lambda_g \lambda_{g-1} = 
(-1)^{g+1} \frac{B_{2g}}{g} \frac{g!}{2^{2g}}
\sum_{p+q=g} \frac1{p! (2p-1)!! \; q! (2q-1)!!}\, .
$$
The right side is,
up to a factor, the sum of even binomial coefficients in the $2g$-th 
line of the Pascal triangle:
$$
\sum_{p+q=g} \frac1{p! (2p-1)!! \; q! (2q-1)!!}
= \sum_{p+q=g} \frac{2^g}{(2p)!\, (2q)!}
= \sum_{p=0}^g \frac{2^g}{(2g)!} {2g \choose 2p} =
\frac{2^{3g-1}}{(2g)!}\, .
$$
Finally, we have
$$
(-1)^{g+1} \frac{B_{2g}}{g} \frac{g!}{2^{2g}} \cdot \frac{2^{3g-1}}{(2g)!} =(-1)^{g+1} \frac{B_{2g}}{2g} \frac1{(2g-1)!!}
$$
as claimed. \qed

\subsection{Local theory of curves}

The local Gromov-Witten theory of curves as developed in \cite{BryPan}
starts with a fundamental pairing with a double ramification cycle
(the resulting matrix then also appears  in the quantum cohomology of
the Hilbert scheme of points of $\C^2$ \cite{OPQ}). We rederive the
double ramification cycle pairing as an almost immediate
consequence of Theorem \ref{FFFF}.

\begin{prop} \label{MMM}
For $g\geq 1$, we have
$$
\int\limits_{\DR_g(a,-a)} \lambda_g \lambda_{g-1} = (-1)^{g+1} \frac{a^{2g}}{(2g)!} \frac{B_{2g}}{2g}.
$$
\end{prop}

Proposition \ref{MMM}
is essentially equivalent to the basic calculation of
\cite[Theorem 6.5]{BryPan} in the Gromov-Witten theory
of local curves. There are two minor differences. 
First, our answer differs from \cite[Theorem 6.5]{BryPan} 
by a factor of  $2g$  coming from the choice of branch point
in \cite[Theorem 6.5]{BryPan}.
Second, if we write out the generating series obtained from
Proposition \ref{MMM} then we only obtain the $\cot(u)$ summand 
in the formula of \cite[Theorem 6.5]{BryPan} after
 a standard correction accounting for the 
geometric difference between the space of rubber maps and the space of maps
relative to three points.

\vspace{10pt}
\noindent{\em Proof.}
We need only compute the contributions to $DR_g(a,-a)$
 of all rational tails type graphs
 since  $\lambda_g \lambda_{g-1}$  vanishes on all other strata. 
Since we just have two marked points, there are exactly two rational tails type graphs:
\begin{enumerate}
\item[$\bullet$]
 the graph with a single vertex,
\item[$\bullet$] the graph with a vertex of genus~$g$, a vertex of genus~0 containing both marked points, and exactly one edge joining the two vertices. 
\end{enumerate}
The contribution of the first graph is (using Lemma~\ref{ggww})
$$
\frac{a^{2g}}{2^g \, g!}\int\limits_{\oM_{g,2}} (\psi_1 + \psi_2)^g \lambda_g \lambda_{g-1} = 
\frac{a^{2g}}{2^g \, g!} \cdot
(-1)^{g+1} \frac{B_{2g}}{2g} \frac1{(2g-1)!!}
=
 (-1)^{g+1} \frac{a^{2g}}{(2g)!} \frac{B_{2g}}{2g}\, .
$$
The contribution of the second graph vanishes because the weight of the edge is equal to 0. 
\qed

For the above computation, we only require Hain's formula for the restriction of the double ramification cycle to $\oM_{g,n}^\ct$ (the extension to $\oM_{g,n}$ is not needed).
In the forthcoming paper \cite{forth}, we will study the Gromov-Witten theory of the rubber over
the resolution of the $A_n$-singularity where the full structure is needed.

\newpage

\begin{center}
\bf{Appendix: \ Polynomiality\  {\em by A. Pixton}} \\
\end{center}

\vspace{10pt}
\noindent{\bf A.1} {\bf Overview.} 
Our goal here is to present a self contained
proof of Proposition $3''$ of Section \ref{Ssec:modr}.
Let $$A=(a_1,\ldots, a_n)$$ be a vector
of $k$-twisted double ramification data for genus $g\geq 0$.
 A longer treatment,
in the context of the deeper polynomiality of $\mathsf{P}_{g}^{d,k}(A)$
in the {\em parts $a_i$ of $A$}, will be given in \cite{PixDR2}.

\vspace{8pt}
\noindent {\bf{Proposition $\bf{3''}$} (Pixton \cite{PixDR2})} {\em 
Let $\Gamma\in\mathsf{G}_{g,n}$ be a fixed stable graph.
Let $\ff$ be a polynomial with $\mathbb{Q}$-coefficients
in variables corresponding to the 
set of half-edges $\mathsf{H}(\Gamma)$. Then the sum 
$$
\mathsf{F}(r) = \sum_{w\in \mathsf{W}_{\Gamma,r,k}} \ff\big(w(h_1), \dots, w(h_{|H(\Gamma)|})\big)
$$
over all possible $k$-weightings mod r 
is (for all sufficiently large $r$)  a {\em polynomial} in $r$ divisible by $r^{h^1(\Gamma)}$.}

\vspace{10pt}
\noindent{\bf A.2}\, {\bf Polynomiality of $\mathsf{F}(r)$.}

\vspace{4pt}
\noindent{\bf A.2.1} {\em Ehrhart polynomials.} The polynomiality of
$\mathsf{F}(r)$ for $r$ sufficiently large will be derived
as a consequence of a small variation of the standard theory \cite{Stan} of 
Ehrhart polynomials of integral convex polytopes.

Let $\mathsf{M}$ be a 
$\v\times \h$ matrix with $\mathbb{Q}$-coefficients satisfying the
following two properties:
\begin{enumerate}
\item[(i)]  $\mathsf{M}$ is {\em totally unimodular}\,:
every square minor of $\mathsf{M}$ has determinant 0, 1, or -1.
\item[(ii)] If  $\mathsf{x} = (x_j)_{1\leq j \leq \h}$ is a column vector
satisfying 
$$\mathsf{Mx} = 0 \ \ \ \text{and}\ \ \   x_j\geq  0\, \ {\text{for all $j$,}} $$
then $\mathsf{x} = 0$. 
\end{enumerate}
We will use the notation $\mathsf{x} \geq 0$ to signify $x_j\geq 0$ for all $j$.
Let
$$\mathbb{Z}^\h_{\geq 0} = \{ \ \mathsf{x}\ | \ \mathsf{x}\in \mathbb{Z}^\h
\  \text{and} \  \mathsf{x}\geq 0\ \}\, .$$

\vspace{8pt}
\noindent {\bf{Proposition $\mathbf{A1}$}} \ {\em 
Let $\mathsf{a},\mathsf{b} \in \mathbb{Z}^\v$ be column vectors with integral coefficients.
Let $\mathsf{Q}\in \mathbb{Q}[x_1,\ldots, x_\h]$ be a polynomial.
Then, the sum over lattice points
\begin{equation*}
\mathsf{S}(r)=\sum_{\mathsf{x}\in \mathbb{Z}^\h_{\geq 0}\, ,\  \mathsf{Mx} = 
\mathsf{a}+r\mathsf{b}} \mathsf{Q}(\mathsf{x})
\end{equation*}
is a {\em polynomial} in the integer parameter $r$ (for all sufficiently 
large $r$).}

\vspace{8pt} 
\noindent {\em Proof.}
We follow the treatment of Ehrhart theory{\footnote{The standard situation for
the Ehrhart polynomial occurs when $\mathsf{Q}=1$ is the
constant polynomial and
$\mathsf{a}= 0$. The variation required here is not very significant, but we
were not able to find an adequate reference.}} in  
\cite[Sections 4.5-4.6]{Stan}.
%

To start, we may assume all coefficients $a_i$ of
$\mathsf{a}$ are nonnegative: if not, multiply the negative
coefficients of $\mathsf{a}$ and the 
corresponding parts of $\mathsf{M}$ and $\mathsf{b}$ by $-1$.

Next, we define a generating function of $\h+\v+1$ variables:
$$\mathsf{R}(X_1, ... , X_\h, Y_1, ... , Y_\v, Z) = 
\sum_{\mathsf{x}\in \mathbb{Z}^\h_{\geq 0}\, ,\  
\mathsf{y}\in \mathbb{Z}^\v_{\geq 0}\, ,\ z \in \mathbb{Z}_{\geq 0}\,, \ 
\mathsf{Mx} = \mathsf{y}+z\mathsf{b}} X^{\mathsf{x}} Y^{\mathsf{y}} Z^z\, .$$
By \cite[Theorem 4.5.11]{Stan}, $\mathsf{R}$  
represents a rational function of the variables $X$, $Y$, and $Z$
 with denominator given by a product of factors of the form 
$$1 - \text{Monomial}(X,Y,Z)\, .$$
The monomials in the denominator factors 
correspond to the extremal rays of the cone in $\mathbb{Z}^{\h+\v+1}_{\geq 0}$
specified by the condition $$\mathsf{Mx} = \mathsf{y}+z\mathsf{b}\, .$$ 
Because $\mathsf{M}$ is totally unimodular, these extremal rays 
all have integral generators with $z$-coordinate equal to 0 or 1. 
Therefore, the denominator is a product of factors of the form 
$$1 - \text{Monomial}(X,Y) \ \ \ \text{and}\ \ \  
1 - \text{Monomial}(X,Y)\cdot Z\, .$$

The sum $\mathsf{S}(r)$ in the statement  of  Proposition {A1} 
 can be obtained from $\mathsf{R}$ by executing the following steps:
\begin{enumerate}
\item[(i)] apply a
differential operator in the $X$ variables to $\mathsf{R}$ which
has the effect of 
 multiplying the coefficient of $X^{\mathsf{x}} Y^{\mathsf{y}} Z^z$ by 
$\mathsf{Q}(\mathsf{x})$, 
\item[(ii)] apply the differential operator 
$$\frac{1}{\prod_{i=1}^{\mathsf{v}}{a_i!}} \,
\frac{\partial^{a_1}}{\partial Y_1^{a_1}}
\frac{\partial^{a_2}}{\partial Y_2^{a_2}} \ldots
\frac{\partial^{a_{\mathsf{v}}}}{\partial Y_{\mathsf{v}}^{a_{\mathsf{v}}}}\,
 ,$$
\item[(iii)] set{\footnote{By property (ii) of the definition of
$\mathsf{M}$, the specialization is well defined.
There are no denominator factors in $\mathsf{R}$ of the form 
$1-\text{Monomial}(X)$.}} $X_j=1$ for all $1\leq j\leq \hH$ and $Y_i=0$ for all $1\leq i \leq \mathsf{v}$,
\item[(iv)] extract the coefficient of $Z^r$.
\end{enumerate}
After step (iii) and before step (iv), 
we have  a rational function in $Z$ with denominator equal to a power of 
$1-Z$. The $Z^r$ coefficient of such a rational function is (eventually) 
a polynomial in $r$. \qed

\vspace{4pt}
\noindent{\bf A.2.1} {\em Vertex-edge matrix.} 

We now prove the polynomiality of
$\mathsf{F}(r)$ for all $r$ sufficiently large.
Let $$A=(a_1,\ldots, a_n)$$ be a vector
of $k$-twisted double ramification data for genus $g\geq 0$.
Let $\Gamma\in\mathsf{G}_{g,n}$ be a fixed stable graph.
For all sufficiently large $r$, the sum over half-edge weightings
\begin{equation}\label{ssxx}
w\in\mathsf{W}_{\Gamma,r,k}\, , \ \ \ 0\leq w(h) < r
\end{equation}
satisfying conditions mod $r$,
\begin{equation}\label{www7}
\mathsf{F}(r) = \sum_{w\in \mathsf{W}_{\Gamma,r,k}} \ff\big(w(h_1), \dots, w(h_{|H(\Gamma)|})\big)\, ,
\end{equation}
can be split into {\em finitely} many sums over 
half-edge weightings satisfying equalities{\footnote{Equalities {\em not}
equalities mod $r$.}} (involving multiples of $r$).

We study each sum in the splitting of \eqref{www7} separately. 
The conditions on the weights $w(h)$ in these sums will be of the following form.
 The inequalities
$$0 \leq w(h) < r$$ 
continue to hold. Moreover,
each leg, edge, and vertex yields an equality:
\begin{enumerate}
\item[(i)] $w(h_i)=a_i+rb_i$ at the $i$th leg $h_i$, where 
$$b_i = 0 \, \text{ if }\, a_i\ge 0 \ \text{ and }\ b_i = 1  \, \text{ if }\,  a_i < 0\, , $$
\item[(ii)] $w(h) + w(h') = 0$\ or\  $w(h) + w(h') = r$\ 
along the edge $(h,h')$,
\item[(iii)] $\sum_{v(h)=v}w(h)= k(2\g(v)-2+\n(v)) + rb_v$\
at the vertex $v$ with $b_v\in \mathbb{Z}$.
\end{enumerate}

The conditions $w(h_i) < r$ are now redundant with the equalities of type (i) and can be removed. Along an edge $(h,h')$, the condition $w(h) < r$ only serves to rule out the possibility $w(h) = r$, $w(h') = 0$ (if we have the equality $w(h) + w(h') = r$). Thus the condition $w(h) < r$ can be removed if we subtract off the contributions with $w(h) = r$, $w(h') = 0$. These contributions are given by more sums of the same form (corresponding to the graph formed by deleting the edge $(h,h')$ from $\Gamma$), multiplied by powers of $r$ and $0$ coming from factors of $w(h)$ and $w(h')$ in the polynomial $\mathsf{Q}$.

Thus the original sum $\mathsf{F}(r)$ can be split into a linear combination (with coefficients being polynomials in $r$) of sums with weight conditions of the following form: 
inequalities
$$w(h) \geq 0$$
and equalities of the form (i), (ii), and (iii) associated to
legs, edges, and vertices as described above.

We must check that
the resulting sums satisfy the hypotheses of Proposition A1.

Consider the graph $\Gamma'$ obtained from $\Gamma$ 
by adding a new vertex in the middle of each edge $e\in \mathsf{E}(\Gamma)$ 
and at the end of each leg $l\in \mathsf{L}(\Gamma)$. 
The half-edge set of $\Gamma$ then 
 corresponds bijectively to the
edge set of $\Gamma'$,
$$\mathsf{H}(\Gamma) \stackrel{\sim}{\longleftrightarrow} \mathsf{E}(\Gamma')\, .$$
We associate {\em each} linear condition on the half-edge weights 
$w(h)$ described in (i-iii) above to the corresponding
{\em vertex} of $\Gamma'$. For types (i) and (ii), the corresponding
vertex of
$\Gamma'$ is a new vertex. For (iii), the corresponding
vertex  already exists in
$\Gamma$.

The matrix $\mathsf{M}$ associated to the linear terms in these
linear conditions is precisely the {\em vertex-edge incidence matrix} of 
$\Gamma'$: $\mathsf{M}$ is a $\v\times \h$ matrix,
$$\v=|\mathsf{V}(\Gamma')|\, , \ \ \h=|\mathsf{E}(\Gamma')| = 
|\mathsf{H}(\Gamma)|\, ,$$
with matrix entry $1$ if the corresponding
vertex of $\Gamma'$ meets the corresponding edge of $\Gamma'$ (and $0$ otherwise).
Since all the entries of $\mathsf{M}$ are nonnegative 
and every column contains a positive entry, we immediately see
 $$\Big(\mathsf{Mx} = 0 \ \ \text{and}\ \ \mathsf{x}\geq 0\Big)
\ \ \Rightarrow  \ \ \mathsf{x} = 0\, .$$
Because $\Gamma'$ is bipartite, the incidence matrix $\mathsf{M}$ is 
totally unimodular by Proposition A2 below.

We may therefore apply Proposition A1 to each sum of the splitting
of \eqref{www7}.
Each such sum is a polynomial in $r$ for sufficiently large $r$.
Hence, $\mathsf{F}(r)$ is a polynomial in $r$ for sufficiently large
$r$. \qed

\vspace{10pt}
\noindent {\bf{Proposition $\mathbf{A2}$}} 
{\em The vertex-edge matrix of a bipartite graph is totally unimodular.}

\vspace{6pt} A proof of general results 
 concerning
 the
total unimodularity of vertex-edge matrices
 can be found in \cite{Hell}. Proposition A2 arises as a special case.

\vspace{10pt}
\noindent{\bf A.3} {\bf Divisibility by $r^{h^1(\Gamma)}$.}
In this section we write
$$ \mathsf{b} = h^1(\Gamma) $$
for convenience. We wish to show that the polynomial $\mathsf{F}(r)$ is divisible by $r^{\mathsf{b}}$. This is implied by checking that the valuation condition,
$$\mathsf{val}_p(\mathsf{F}(p)) \geq \mathsf{b}\, ,$$
holds for {\em all} sufficiently large primes $p$. Here,
$$\mathsf{val}_p: \mathbb{Q} \rightarrow \mathbb{Z}$$
is the $p$-adic valuation. 

Let $\hH=|\mathsf{H}(\Gamma)|$ be the number of half-edges of $\Gamma$.
We write the
 evaluation $\mathsf{F}(p)$ as
$$\mathsf{F}(p) = \sum_{0 \leq  w_1, \ldots , w_\hH \leq p-1} 
\ff(w_1,\ldots , w_\hH)\, $$
where the sum is over  $\hH$-tuples $(w_1, \ldots , w_\hH)$ 
subject to $\hH-\mathsf{b}$ linear conditions mod $p$. The polynomial
$\ff$ has  ($p$-integral) $\mathbb{Q}$-coefficients. 

We count the linear conditions mod $p$ as follows. By the definition of a $k$-weighting mod $p$, we start with
$$\mathsf{v}=|\mathsf{L}(\Gamma)|+|\mathsf{E}(\Gamma)|+|\mathsf{V}(\Gamma)|\ $$
linear conditions, but the conditions are not independent since 
the vector $A$ of $k$-twisted double ramification data satisfies
$$\sum_{i=1}^n a_i = k(2g-2+n)\, .$$
After removing a redundant linear relation, we have $\mathsf{v}-1$ independent linear conditions.
 Since
$$\hH= |\mathsf{H}(\Gamma)| = |\mathsf{L}(\Gamma)|+2|\mathsf{E}(\Gamma)|\ ,$$
and $\Gamma$ is connected, we see
$$\hH-\mathsf{b} = \hH-(\mathsf{E}(\Gamma)-\mathsf{V}(\Gamma)+1) = \mathsf{v}-1\, .$$

Let $\fff$ be a polynomial with $\mathbb{Q}$-coefficients
in $\hH$ variables, with all coefficients $p$-integral.
We will prove  {\em every}  summation 
\begin{equation}\label{n444}
\widehat{\mathsf{F}}(p) = 
\sum_{0 \leq  w_1, \ldots , w_\hH \leq p-1} \fff(w_1,\ldots , w_\hH)\, ,
\end{equation}
where the sum is over  $\hH$-tuples $(w_1, \ldots , w_\hH)$ 
subject to $\hH-\mathsf{b}$ linear conditions mod $p$,
is divisible by $p^{\mathsf{b}}$ for $p \geq \text{deg}(\fff)+2$.

We may certainly take $\fff$ to be a monomial. 
Without loss of generality{\footnote{To increase $w_1$ to $w_1^2$, add a new variable $w_{\hH+1}$ and
a new linear condition $w_1=w_{\hH+1}$.
If a variable $w_j$ does not appear in $\widehat{\mathsf{Q}}$, then a reduction to a
fewer variable case can be made.}}}, we
need only consider
$$\fff(w_1, \ldots , w_\hH) = \prod_{j=1}^\hH w_j\, .$$ 
We write the linear conditions
imposed
in the sum \eqref{n444} as
$$\sum_{j=1}^\hH \bb_{ij} w_j + \cc_i = 0 \mod p\, $$
for $1\leq i \leq \hH-\mathsf{b}$.
Alternatively, we may impose the conditions by standard $p$-th root of unity 
filters:
\begin{equation*}p^{\hH-\mathsf{b}}\, \widehat{\mathsf{F}}(p) = 
\sum_{z_1, \ldots , z_{\hH-\mathsf{b}}}\,
\sum_{0 \leq w_1, \ldots , w_\hH \leq  p-1}\,  
\prod_{j=1}^\hH w_j
\cdot
\prod_{\substack{1 \leq i \leq \hH-\mathsf{b} \\ 1 \leq j \leq \hH}} z_i^{\bb_{ij} w_j} \cdot
\prod_{1 \leq i \leq \hH-\mathsf{b}} 
z_i^{\cc_i}\, .
\end{equation*}
The first sum on the right is over all $(\hH-\mathsf{b})$-tuples
$({z_1, \ldots , z_{\hH-\mathsf{b}}})$
of $p$-th roots of unity. 
The second sum on the right is over {\em all} $\hH$-tuples
$(w_1,\ldots,w_\hH)$.

Consider the following  polynomial with integer coefficients:
$$G(Z) = \sum_{0 \leq  a \leq  p-1} aZ^a = 
\frac{Z - pZ^p + (p-1)Z^{p+1}}{(1 - Z)^2}\, .$$
We can use $G$ to rewrite the 
terms inside the first sum in the above expression:
\begin{equation}\label{w999}
p^{\hH-\mathsf{b}}\, \widehat{\mathsf{F}}(p)
 = 
\sum_{z_1, \ldots , z_{\hH-\mathsf{b}}}\, 
 \prod_{1 \leq j \leq \hH} 
G\left(\prod_{1 \leq i \leq \hH-\mathsf{b}} z_i^{\bb_{ij}}\right) \cdot
\prod_{1 \leq i \leq \hH-\mathsf{b}} z_i^{\cc_i}\, .
\end{equation}

For every $p$-th root of unity $z$, $G(z)$ 
has $p$-adic valuation{\footnote{We view the $p$-th root of unity $z$
as lying in $\overline{\mathbb{Q}}_p$. The $p$-adic valuation of $x=z-1$
can be calculated from the minimal polynomial
$$\sum_{i=0}^{p-1} (x+1)^i = 0\ .$$}}
at least $\frac{p-2}{p-1}$:
either $z=1$ and $G(z)= \frac{p(p-1)}{2}$ or
$$z\neq 1 \ \ \Rightarrow \ \  G(z)= \frac{p}{z-1} \ \ \text{and}\ \ 
\mathsf{val}_p(z-1)=\frac{1}{p-1}\, .$$
Therefore, every term in the sum on the right side of \eqref{w999}
 has $p$-adic valuation at least $\hH\cdot \frac{p-2}{p-1}$. Hence,
$$\mathsf{val}_p(\widehat{\mathsf{F}}(p)) \geq  \hH\cdot \frac{p-2}{p-1} 
- (\hH-\mathsf{b})\, .$$
Since $\mathsf{val}_p(\widehat{\mathsf{F}}(p))$ is an {\em integer},
$$\mathsf{val}_p(\widehat{\mathsf{F}}(p)) \geq \mathsf{b}$$
 if $p \geq \hH+2$. \qed

\vspace{10pt}
\noindent{\bf A.4} {\bf Polynomiality of $\mathsf{P}_g^{d,k}(A)$ in the $a_i$.}
For double ramification data of fixed length $n$, consider
$$\mathsf{P}_g^{d,k}(a_1,\ldots,a_n) \in R^d(\oM_{g,n})$$
as a function of 
$$(a_1,\ldots,a_n)\in \mathbb{Z}^n\, , \ \ \ \sum_{i=1}^{n} a_i =k(2g-2+n)\, $$
holding $g$, $d$, and $k$ fixed. 
The following property is proven in
\cite{PixDR2}:  

\vspace{3pt}
\noindent \hspace{60pt} $\mathsf{P}_g^{d,k}(a_1,\ldots,a_n)$ {\em is polynomial
in the variables $a_i$}.

\noindent By Theorem \ref{FFFF}, the double ramification cycle $\mathsf{DR}_g(A)$
is then also polynomial in the variables $a_i$. 
Shadows of the polynomiality of the double ramification cycle
were proven earlier in \cite{BSSZ,OP1}.

\newpage
\newcommand{\arXiv}[1]{\texttt{arXiv:#1}}

\vspace{+12 pt}

\noindent Institut Math\'ematique de Jussieu\\
\noindent felix.janda@imj-prg.fr  

\vspace{+6pt}
\noindent Departement Mathematik, ETH Z\"urich \\
\noindent rahul@math.ethz.ch

\vspace{+6 pt}
\noindent
Department of Mathematics, MIT\\
apixton@mit.edu

\vspace{+6 pt}
\noindent
CNRS, Institut Math\'ematique de Jussieu\\
zvonkine@math.jussieu.fr


\begin{thebibliography}{99}

\bibitem{AGV} D. Abramovich, T. Graber, and A. Vistoli, {\em Gromov-{W}itten theory of Deligne-Mumford stacks}, Amer. J. Math. {\bf 130} (2008), 1337--1398. 

\bibitem{BryPan} J. Bryan and R. Pandharipande, {\em The local Gromov-Witten theory of curves}, JAMS {\bf 21} (2008), 101--136. 

\bibitem{BSSZ} A. Buryak, S. Shadrin, L. Spitz, and D. Zvonkine, {\em Integrals of $\psi$-classes over double ramification cycles}, Amer. J. Math. {\bf 137} 
(2015), 699--737. 

\bibitem{Cav} R. Cavalieri, S. Marcus, and J. Wise, {\em Polynomial
families of tautological classes on $\cM_{g,n}^{\mathsf{rt}}$}, J. Pure and
Applied Algebra {\bf 216} (2012), 950--981.

\bibitem{Chiodo} A. Chiodo, {\em Stable twisted curves and their $r$-spin structures}, Ann. Inst. Fourier  {\bf 58} (2008), 1635--1689. 

\bibitem{Chiodo2} A. Chiodo, {\em Towards an enumerative geometry of the moduli space of twisted curves and rth roots}, Compos. Math. {\bf 144} (2008), 1461--1496. 

\bibitem{ChiZvo} A. Chiodo and D. Zvonkine, {\em Twisted Gromov--Witten $r$-spin potential and Givental's quantization}, 
Advances in Mathematical Physics {\bf 13} (2009), 1335-1369.


\bibitem{cj} E. Clader and F. Janda, {\em Pixton's double ramification cycle relations}, arxiv:1601.02871.

\bibitem{Fab} C. Faber, {\em A conjectural description of the tautological ring of the moduli space of curves} in {\em Moduli of curves and abelian varieties}, 109--129, Aspects Math., E33, Vieweg, Braunschweig, 1999. 

\bibitem{FP-H} C. Faber and R. Pandharipande, {\em Hodge integrals
and Gromov-Witten theory}, Invent. Math. {\bf 139} (2000), 173--199.

\bibitem{FP} C. Faber and R. Pandharipande, {\em Relative maps and tautological classes}, JEMS {\bf 7} (2005), 13--49. 

\bibitem{FP-Handbook} C. Faber  and R. Pandharipande, {\em Tautological
and non-tautologi\-cal cohomology of the moduli
space of curves}
in {\em Handbook of moduli. Vol.~I}, 293--330, Adv. Lect. Math. (ALM), {\bf 24}, Int. Press, Somerville, MA, 2013. 

\bibitem{FarP} G. Farkas and R. Pandharipande, {\em The moduli space
of twisted canonical divisors} with Appendix by
F. Janda, R. Pandharipande, A. Pixton, and D. Zvonkine, \arXiv{1508.07940}.

\bibitem{Geer} G. van der Geer, in {\em Cycles on the
moduli space of  abelian varieties}, 65--89, Aspects Math., E33, Vieweg, Braunschweig, 1999. 

\bibitem{GetP} E. Getzler and R. Pandharipande, {\em Virasoro constraints and the Chern classes of the Hodge bundle}, Nuclear Phys. B {\bf 530} (1998), 701--714. 

\bibitem{GrP} T. Graber and R. Pandharipande, {\em Localization of virtual classes}, Invent. Math. {\bf 135} (1999), 487--518. 

\bibitem{GrP2} T. Graber and R. Pandharipande, {\em Constructions of nontautological classes on moduli spaces of curves}, Michigan Math. J. {\bf 51} (2003), 93--109. 

\bibitem{GrV} T. Graber and R. Vakil, {\em Relative virtual localization and vanishing of tautological classes on moduli spaces of curves}, Duke Math. J. {\bf 130} (2005), 1--37. 


\bibitem{GruZak} S. Grushevsky and D. Zakharov, {\em The double ramification cycle and the theta divisor}, Proc. Amer. Math. Soc. {\bf 142} (2014), 4053--4064. 

\bibitem{Hain} R. Hain, {\em Normal functions and the geometry of moduli spaces of curves} in {\em Handbook of moduli. Vol.~I}, 527--578, Adv. Lect. Math. (ALM), {\bf 24}, Int. Press, Somerville, MA, 2013. 

\bibitem{Hell} I. Heller and C. Tompkins, {\em An extension of a theorem 
of Dantzig's}, in  
{\em Linear inequalities and related systems} by
H. Kuhn and A. Tucker, Annals of Mathematics Studies 38,  Princeton University Press (1956),  247--254.


\bibitem{J} F. Janda, {\em Relations on $\oM_{g,n}$ via equivariant Gromov-Witten theory of $\PP^1$}, \arXiv{1509.08421}.

\bibitem{forth} F. Janda, R. Pandharipande, A. Pixton, D. Zvonkine, {\em in preparation}.

\bibitem{Jarvis} T. J. Jarvis, {\em Geometry of the moduli of higher spin curves,} Int. J. Math. {\bf 11} (2000), 637--663. 

\bibitem{JPT} P. Johnson, R. Pandharipande, and  H.-H. Tseng, {\em Abelian Hurwitz-Hodge integrals}, {\em Abelian Hurwitz-Hodge integrals}, Michigan Math. J.
{\bf 60} (2011), 171--198.
\arXiv{0803.0499}

\bibitem{JLi} J.~Li, {\em Lecture notes on relative GW-invariants,} \\
\small{\texttt{http:\slash\slash{}users.ictp.it\slash$\sim\!\!$
pub\_off\slash{}lectures\slash{}lns019\slash{}Jun\_Li\slash{}Jun\_Li.pdf}}

\bibitem{Li2} J.~Li, {\em A degeneration formula of GW invariants,} 
J.~Differential Geom. {\bf 60} (2002), 199--293.

\bibitem{MW} S. Marcus and J. Wise, {\em Stable maps to
rational curves and the relative Jacobian}, \arXiv{1310.5981}.


\bibitem{MauPan} D.~Maulik and R.~Pandharipande, {\em A topological view of Gromov-Witten theory,} Topology {\bf 45} (2006), 887--918.

\bibitem{Mum} D. Mumford, {\em Towards an enumerative geometry of
the moduli space of curves}, in {\em Arithmetic and Geometry}
(M. Artin and J. Tate, eds.), Part II, Birkh\"auser, 1983, 271--328.

\bibitem{OP1} A. Okounkov and R. Pandharipande, {\em Gromov-Witten
theory, Hurwitz numbers, and completed cycles}, Ann. of Math. 
{\bf 163} (2006), 517--560. 

\bibitem{OPQ} A. Okounkov and R. Pandharipande, {\em Quantum 
cohomology of the Hilbert scheme of points of the plane}, Invent. Math. 
{\bf 179} (2010), 523--557. 

\bibitem{PP} R. Pandharipande and A. Pixton, {\em Relations in the
tautological ring of the moduli space of curves}, \arXiv{1301.4561}.

\bibitem{PPZ} R. Pandharipande, A. Pixton, and D. Zvonkine, {\em Relations on $\oM_{g,n}$ via 3-spin structures}, JAMS {\bf 28} (2015), 279--309.

\bibitem{P} A.~Pixton, {\em Conjectural relations in the tautological ring of $\oM_{g,n}$}, \\ \arXiv{1207.1918}.

\bibitem{PixDR} A.~Pixton, {\em Double ramification cycles and tautological relations on $\oM_{g,n}$}, available from the author.

\bibitem{PixDR2} A. Pixton, {\em On combinatorial properties of the explicit expression for double ramification cycles}, in preparation.

\bibitem{R} B. Riemann, {\em Theorie der Abel'schen Functionen}, J. Reine Angew. Math. {\bf 54} (1857), 115--155. 

\bibitem{Stan} R. Stanley, {\em Enumerative combinatorics: Vol I}, Cambridge
Univ. Press: Cambridge, 1999.
\end{thebibliography}
\end{document}